\documentstyle[amstex,amscd]{article}
\title{Extension of 2-forms and symplectic varieties}
\author{Yoshinori Namikawa}
\date{ }
\begin{document}
\maketitle
\begin{center}
{\bf Introduction}
\end{center}

In this paper we shall prove two theorems 
(Stability Theorem, Local Torelli Theorem) 
for symplectic varieties.
  
 Let us recall the notion of a symplectic singularity. 
Let $X$ be a good representative 
of a normal singularity. 
Then the singularity is symplectic 
if the regular locus $U$ of $X$ admits 
an everywhere non-degenerate holomorphic closed 2-form 
$\omega$ where $\omega$ extends 
to a regular form on $Y$ for a resolution of 
singularities $Y \to X$.  
Similarly we say that a normal compact Kaehler space 
$Z$ is a symplectic variety 
if the regular locus $V$ of $Z$ admits a 
non-degenerate holomorphic closed 2-form 
$\omega$ where $\omega$ extends to a regular 
form on $\tilde Z$, 
where $\tilde Z \to Z$ is a resolution of 
singularities of $Z$. 
When $Z$ has a resolution $\pi: \tilde Z \to Z$ 
such that $(\tilde Z, \pi^*\omega)$ is a 
symplectic manifold, 
we call $Z$ has a symplectic resolution.  
\vspace{0.15cm}

{\bf Examples} 

(i) This is one of examples of symplectic singularities 
studied in [Be 1]. For details see [Be 1] 
and the references there. 
Let $Q \subset {\bold P}^{n-1}$ be a general 
quadratic hypersurface. 
Identify a point of the Grassmannian $Gr(2, n)$ 
with a line in ${\bold P}^{n-1}$. 
Let $Gr_{iso}(2, n)$ be the subvariety of $Gr(2, n)$ 
corresponding to the lines of ${\bold P}^{n-1}$ 
contained in $Q$. 
It is checked that 
$\dim Gr_{iso}(2, n) = \dim Gr(2, n) - 3 = 2n - 7$. 
Embed $Gr_{iso}(2, n)$ into ${\bold P}^{1/2n(n-1)-1}$ 
by the Pl\"ucker embedding $Gr(2, n) \to 
{\bold P}^{1/2n(n-1)-1}$. 
Now consider the cone $X$ over $Gr_{iso}(2, n)$. 
Then the germ $(X, 0)$ at the vertex is a 
symplectic singularity
 of dimension $2n - 6$. 
The $X$ is actually obtained as the closure 
$\bar{\cal O}_{min}$ of the minimal nilpotent orbit 
${\cal O}_{min}$ of the Lie algebra $Lie(SO(n))$, 
and ${\cal O}_{min}$ has the Kostant-Kirillov 
symplectic 2-form. 
\vspace{0.12cm}

(ii) Let $A := {\bold C}^{2l}/\Gamma$ be an 
Abelian variety 
of dimension $2l$. 
Let $(z_1, z_2, ..., z_{2l-1}, z_{2l})$ be the 
standard coordinates 
of ${\bold C}^{2l}$. Then ${\bold Z}/2{\bold Z}$ acts on 
$A$ by $z_i \to -z_i$ $(i = 1, ..., 2l)$. 
The quotient $Z$ of $A$ by the action becomes a 
symplectic variety 
of dimension $2l$. A symplectic 2-form is, for example, 
given by 
$\Sigma_{1 \leq i \leq l}dz_i \wedge dz_{l+i}$. 
The $Z$ has singularities, 
and $Z$ has no symplectic resolution when $l > 1$. 
\vspace{0.12cm}

(iii) These are symplectic varieties 
studied by O'Grady [O]. 
Let $S$ be a polarized K3 surface. 
Let $c$ be an even number with $c \geq 4$. 
Denote by $\overline{M}_{0,c}$ the moduli space of 
rank 2 semi-stable torsion free sheaves with 
$c_1 = 0$ and $c_2 = c$. 
$\overline{M}_{0,c}$ becomes a projective 
symplectic variety 
of $\dim = 4c - 6$. 
The singular locus $\Sigma$ has dimension $2c$. 
Moreover, O'Grady showed that $\overline{M}_{0,4}$ 
has a symplectic resolution, 
however $\overline{M}_{0,c}$ has ${\bold Q}$-factorial 
terminal singularities when $c \geq 6$ (cf. section 3 
of the e-print version of [O]: alg-geom/9708009).   
Therefore $\overline{M}_{0,c}$ have no 
symplectic resolution 
when $c \geq 6$.  \vspace{0.12cm}

 A symplectic singularity / variety will play 
an important role in 
the generalized Bogomolov decomposition conjecture 
(cf. [Kata], [Mo]): \vspace{0.12cm}

{\bf Conjecture}: {\em Let $Y$ be a smooth 
projective variety 
over ${\bold C}$ with Kodaira dimension 0. 
Then there is a finite etale cover $Y' \to Y$ 
such that $Y'$ is birationally equivalent to 
$Y_1 \times Y_2 \times Y_3$, where $Y_1$ 
is an Abelian variety, $Y_2$ is a symplectic variety, 
and $Y_3$ is a Calabi-Yau variety. }  \vspace{0.12cm}

In this conjecture we hope that 
it is possible to replace 
$Y_2$ and $Y_3$ by their birational models 
with only ${\bold Q}$-factorial 
terminal singularities respectively. 
Main results are these. \vspace{0.15cm}
 
{\bf Theorem 7}({\bf Stability Theorem}): 
{\em Let $(Z, \omega)$ be a projective 
symplectic variety.  
Let $g : {\cal Z} \to \Delta$ be a 
projective flat morphism 
from ${\cal Z}$ to a 1-dimensional unit disc 
$\Delta$ with $g^{-1}(0) = Z$. 
Then $\omega$ extends sideways 
in the flat family so that it gives a 
symplectic 2-form 
$\omega_t$ on each fiber $Z_t$ for 
$t \in \Delta_{\epsilon}$ 
with a sufficiently small $\epsilon$.} 
\vspace{0.15cm}

In the above, the result should also hold for a 
(non-projective) symplectic variety $(Z, \omega)$ 
and for a proper flat morphism $g$. 
But two ingredients remained unproved in 
the general case 
(cf. Remark below Theorem 7). \vspace{0.12cm}

Let $Z$ be a symplectic variety.  
Put $\Sigma := \mathrm{Sing}(Z)$ and $U := Z 
\setminus \Sigma$.  
Let $\overline \pi : {\cal Z} \to S$ 
be the Kuranishi family of $Z$, which is, by definition, 
a semi-universal flat deformation of $Z$ with 
$\overline \pi^{-1}(0) = Z$ for the reference point 
$0 \in S$. 
When $\mathrm{codim}(\Sigma \subset Z) \geq 4$, 
$S$ is smooth by [Na 1, Theorem 2.4]. 
${\cal Z}$ is not projective over $S$. 
But we can show that every member 
of the Kuranishi family is a symplectic 
variety (cf. Theorem 7'). 
Define ${\cal U}$ to be the locus in ${\cal Z}$ 
where $\overline \pi$ is a smooth map and 
let $\pi : {\cal U} \to S$ be the restriction of 
$\overline \pi$ to ${\cal U}$. Then we have \vspace{0.2cm}

{\bf Theorem 8}({\bf Local Torelli Theorem}): 
{\em Assume that $Z$ is a ${\bold Q}$-factorial 
projective symplectic variety. Assume 
$h^1(Z, {\cal O}_Z) = 0$, $h^0(U, \Omega^2_U) = 1$, 
$\dim Z = 2l \geq 4$ and 
$\mathrm{Codim}(\Sigma \subset Z) \geq 4$. 
Then the following hold.} \vspace{0.15cm}

(1) {\em $R^2\pi_*(\pi^{-1}{\cal O}_S)$ 
is a free ${\cal O}_S$ module of finite rank. 
Let ${\cal H}$ be the image of the composite 
$R^2{\overline \pi}_*{\bold C} 
\to R^2\pi_*{\bold C} \to 
R^2\pi_*(\pi^{-1}{\cal O}_S)$. 
Then ${\cal H}$ is a local system on $S$ 
with ${\cal H}_s = H^2({\cal U}_s, {\bold C})$ 
for $s \in S$.} \vspace{0.15cm}

(2) {\em  The restriction map 
$H^2(Z, {\bold C}) \to H^2(U, {\bold C})$ 
is an isomorphism. 
Take a resolution $\nu: \tilde Z \to Z$ 
in such a way that 
$\nu^{-1}(U) \cong U$. For $\alpha 
\in H^2(U, {\bold C})$ 
we take a lift $\tilde\alpha \in 
H^2(\tilde Z, {\bold C})$ 
by the composite $H^2(U, {\bold C}) 
\cong H^2(Z, {\bold C}) 
\to H^2(\tilde Z, {\bold C})$. 
Choose $\omega \in H^0(U, \Omega^2_U) = 
{\bold C}$. This $\omega$ extends to a 
holomorphic 2-form 
on $\tilde Z$. Normalize $\omega$ 
in such a way that 
${\int_{\tilde Z} (\omega \overline \omega)^l} = 1$. 
Then one can define a quadratic form 
$q : H^2(U, {\bold C}) \to {\bold C}$ as} 

$$  q(\alpha) :=  
l/2 {\int_{\tilde Z}(\omega \overline \omega)^{l-1}
\tilde\alpha^2} 
+ (1-l)(\int_{\tilde Z}\omega^l\overline\omega^{l-1}
\tilde\alpha)
(\int_{\tilde Z}\omega^{l-1}
\overline\omega^l\tilde\alpha ). $$

{\em This form is independent of the choice of 
$\nu : \tilde Z \to Z$.} \vspace{0.15cm}
 
(3) {\em Put $H := H^2(U, {\bold C})$. 
Then there exists a trivialization of the local system 
${\cal H}$:  ${\cal H} \cong H \times S$. 
Let ${\cal D} := \{ x \in {\bold P} (H) ; 
q(x) =0, q(x + \overline x) > 0 \}$. 
Then one has a period map $p: S \to {\cal D}$ and $p$ 
is a local isomorphism.} \vspace{0.2cm}

Note that when $Z$ is a symplectic variety 
with terminal singularities, 
the condition $\mathrm{codim}(\Sigma \subset Z) \geq 4$ 
is always satisfied ([Na 2]).
 
The stability theorem will be proved 
by using the following theorem and 
the fact that a projective variety 
with rational singularities has Du Bois 
singularities (cf. [Ko]):  \vspace{0.15cm}

{\bf Theorem 4.}. {\em  Let $X$ be a Stein open subset 
of a complex algebraic variety. 
Assume that $X$ has only rational 
Gorenstein singularities. 
Let $\Sigma$ be the singular locus of $X$ 
and let $f : Y \to X$ be a resolution of singularities 
such that $f\vert_{Y\setminus f^{-1}(\Sigma)} 
: Y \setminus f^{-1}(\Sigma) \cong X \setminus \Sigma$. 
Then $f_*\Omega^2_Y \cong i_*\Omega^2_U$ 
where $U := X \setminus \Sigma$ and 
$i : U \to X$ is a natural injection.} \vspace{0.15cm}

The same result was obtained by 
D. van Straten and Steenbrink [S-S] 
for an arbitrary isolated normal singularity 
with $\dim \geq 4$, and later Flenner [Fl] proved it 
for arbitrary normal singularity 
with $\mathrm{Codim}(\Sigma \subset X) \geq 4$. 
By Theorem 4 $X$ is a symplectic singularity 
if and only if $X$ is rational Gorenstein 
and the regular part of $X$ 
admits an everywhere non-degenerate 2-form 
(cf. Theorem 6). 
This fact is often useful to determine 
that certain kinds of singularities given 
by G.I.T. quotient are symplectic (Example 6').    

The local Torelli theorem has been proved 
for non-singular symplectic varieties  
by Beauville [Be 2, Theoreme 5]. 
In our singular case, it is based on the 
Hodge decomposition 
$H^2(U, {\bold C}) = H^0(U, \Omega^2_U) 
\oplus H^1(U, \Omega^1_U) \oplus H^2(U, {\cal O}_U)$. 
We need the condition that 
$\mathrm{codim}(\Sigma \subset Z) \geq 4$ 
to have this decomposition (cf. [Oh, Na 1, Lemma 2.5]). 
We shall prove that $R^2\overline\pi_*{\bold C}$ 
is a constant sheaf around $0 \in S$ by using the 
${\bold Q}$-factoriality of $Z$. 
When $Z$ is not ${\bold Q}$-factorial, 
the statement for ${\cal H}$ in (1) does not hold 
as it stands because $H^2(Z, {\bold C}) 
\to H^2(U, {\bold C})$ 
is not surjective. 

We have formulated the local Torelli theorem 
for a projective symplectic variety, but it 
is possible to get similar statements for 
a general (non-projective) symplectic variety. 
For a non-projective symplectic variety, 
${\bold Q}$-factoriality should be replaced by 
a certain condition which is equivalent to 
the ${\bold Q}$-factoriality in the 
projective case (cf. Remark (2) on the 
final page). 

The following problem would be 
of interest in view of Global Torelli Problem. 
\vspace{0.15cm}

{\bf Problem}. {\em Let $Z$ be a ${\bold Q}$-factorial 
projective symplectic variety with terminal 
singularities. 
Assume that $Z$ is smoothable by a 
suitable flat deformation. Is $Z$ then non-singular? } 
\vspace{0.15cm}

In the first section we shall prove Theorem 4. 
Other two theorems are proved in 
the second section. \vspace{0.15cm}   

{\bf Notation}.  Let ${\cal F}$ be a coherent sheaf 
on a normal crossing variety $D$. 
Assume that ${\cal F}$ is a locally free sheaf 
on the regular locus of $D$. 
Then $\hat{\cal F} = {\cal F}/(torsion)$ by definition. 
Here $(torsion)$ means the subsheaf of the sections 
whose support are contained in 
the singular locus of $D$.   
\vspace{0.15cm}

\begin{center}
{\bf \S 1. Extension properties for 
Rational Gorenstein Singularities}
\end{center}
\vspace{0.15cm}

{\bf Proposition 1}. {\em Let $X$ be a Stein open subset 
of a complex algebraic variety.  
Assume that $X$ has only rational 
Gorenstein singularities. 
Let $\Sigma$ be the singular locus of $X$ and 
let $f : Y \to X$ be a resolution of singularities 
such that $f\vert_{Y\setminus f^{-1}(\Sigma)} : 
Y \setminus f^{-1}(\Sigma) \cong X \setminus \Sigma$ 
and $D := f^{-1}(\Sigma)$ is a simple normal crossing 
divisor. 
Then $f_*\Omega^2_Y(\log D) \cong i_*\Omega^2_U$ 
where $U := X \setminus \Sigma$ and $i : U \to X$ 
is a natural injection.}  \vspace{0.15cm}

{\em Proof}. Let $\omega \in H^0(U, \Omega^2_U)$. 
$\omega$ 
is a meromorphic 2-form on $Y$. 
In fact, $\mathrm{Coker}[f_*\Omega^2 
\to i_*\Omega^2_U]$ 
is a torsion sheaf whose support is contained in 
$\Sigma$. 
Hence $\phi\omega$ is an element of 
$\Gamma(X, f_*\Omega^2_Y)$ for a suitable holomorphic 
function $\phi$ on $X$.  

Since $(X,p) \cong (R.D.P) \times ({\bf C}^{n-2}, 0)$ 
for all $p \in X$ outside certain codimension 3 (in $X$) 
locus $\Sigma_0 \subset \Sigma$ ([Re]), 
it is clear that $f_*\Omega^2_Y(\log D) 
\cong i_*\Omega^2_U$ 
at such points $p$. Let $F$ be an irreducible component 
of $D$ with $f(F) \subset \Sigma_0$. 
Put $k := \dim \Sigma_0 - \dim f(F)$. 
We shall prove that $\omega$ has at worst 
log pole along $F$ by the induction on $k$.  

(a) $k = 0$:

(a-1):  Put $l := \mathrm{codim} (\Sigma_0 \subset X)$. 
Note that $l \geq 3$. Take a general $l$ dimensional 
complete intersection $H := H_1 \cap ... \cap H_{n-l}$. 
Let $p \in H \cap f(F)$. Since $H$ is general, 
$p \in f(F)$ is a smooth point. 
Replace $X$ by a suitable small open neighborhood of $p$. 
Then $H \cap f(F) = \{p\}$. 
$H$ has a unique dissident point $p$ and other 
singularities 
are locally isomorphic to $(R.D.P.)
\times({\bold C}^{l-2},0)$. 
By perturbing $H$ we can define a flat holomorphic map 
$g : X \to \Delta^{n-l}$ such that the fiber $X_0$ 
over $0 \in \Delta^{n-l}$ coincides with $H$. 
We may assume that $g$ has a section passing through 
$p$ and each fiber $g^{-1}(t)$ intersects $f(F)$ 
only in this section. 
The map $f : Y \to X$ gives a simultaneous resolution of 
$X_t$ ($t \in \Delta^{n-l}$). 
Since $H$ is general and $X$ is sufficiently small, 
$D_t := D \cap Y_t$ are normal crossing divisors of $Y_t$ 
for all $t \in \Delta^{n-l}$. Let $D'$ 
be the union of irreducible components of 
$D$ which are mapped in this section. 
$D' \to \Delta^{n-l}$ is a proper map. 
We put $\pi = g \circ f$. 
We often write $\Delta$ for $\Delta^{n-l}$.

(a-2):  We shall prove the following. \vspace{0.15cm}

{\bf Claim}. {\em By replacing $\Delta^{n-l}$ 
by a smaller disc and by restricting everything 
(e.g. $X$, $Y$, $D$, $D'$ ...) over the new disc, 
we have a subset $K \subset Y$ which contains $D'$ 
and which is proper over $\Delta^{n-l}$ 
with the following property:}  

{\em The $\omega$ is mapped to zero 
by the comoposition of the maps}

$$   H^0(U, \Omega^2_U) \cong 
H^0(Y \setminus D', \Omega^2_{Y \setminus D'}(\log D)) 
\to H^0(Y \setminus K, \Omega^2_{Y \setminus K}(\log D)) 
\to H^1_K(Y, \Omega^2_Y(\log D)). $$    \vspace{0.2cm}

If the claim is verified, then 
$\omega\vert_{Y \setminus K}$ extends to a 
logarithmic 2-form on $Y$. It is clear that 
its restriction to $U$ is $\omega$.  

{\em Proof of Claim}.  We shall first prove that 
$R^1\pi_{!}\Omega^2_Y(\log D) = 0$. 
There is a filtration 
$\pi^*\Omega^2_{\Delta} \subset 
{\cal G} \subset \Omega^2_Y(\log D)$ 
which yields exact sequences 

$$   0 \to {\cal G} \to \Omega^2_Y(\log D) 
\to \Omega^2_{Y/\Delta}(\log 
D) \to 0 $$

$$   0 \to \pi^*\Omega^2_{\Delta} \to {\cal G} \to 
\pi^*\Omega^1_{\Delta}\otimes 
\Omega^1_{Y/\Delta}(\log D) \to 0. $$

By the exact sequences it suffices to prove 
that $R^1\pi_{!}{\cal O}_Y 
= R^1\pi_{!}\Omega^1_{Y/\Delta}(\log D) = 
R^1\pi_{!}\Omega^2_{Y/\Delta}(\log D) = 0$.  

We shall use the relative duality theorem due to 
Ramis and Ruget [R-R] 
to prove these facts. Before applying 
the relative duality we note 
that $R^i\pi_*\Omega^{l-2}_{Y/\Delta}(\log D)(-D) = 
R^i\pi_*\Omega^{l-1}_{Y/\Delta}(\log D)(-D) = 
R^i\pi_*\Omega^l_{Y/\Delta}(\log D)(-D) = 0 $ for 
$i \geq l - 1$. To prove these, 
we only have to check that 
$R^i\pi_*\Omega^{l-2}_{Y_t}(\log D_t)(-D_t) = 
R^i\pi_*\Omega^{l-1}_{Y_t}(\log D_t)(-D_t) = 
R^i\pi_*\omega_{Y_t} = 0$ 
for $i \geq l -1$ and for $t \in \Delta$, 
for example, by using [B-S, VI, Cor. 4.5, (i)]).  
Since $l \geq 3$, these follow from a vanishing 
theorem in [St] except for the vanishing of 
$R^2\pi_*\Omega^1_{Y_t}(\log D_t)(-D_t)$. 
But, by the same argument as 
the proof of [Na-St, Theorem (1.1)] we see that 
$R^2\pi_*\Omega^1_{Y_t}(\log D_t)(-D_t) = 0$. 

Now the relative duality says that 

$$ {\bold R}\pi_{!}{\bold R}
{\cal H}om(Y; \Omega^{l-j}_{Y/\Delta}(\log D)(-D)
\otimes\pi^*\Omega^{n-l}_{\Delta}, \omega^{\cdot}_Y) $$
 
$$\cong {\bold R}{\cal H}omtop
(\Delta, 
{\bold R}\pi_*(\Omega^{l-j}_{Y/\Delta}(\log D)(-D)
\otimes\pi^*\Omega^{n-l}_{\Delta}), 
\omega^{\cdot}_{\Delta}) $$
for $j = 0, 1, 2$. \vspace{0.15cm}
          
We have $R^{1-n}\pi_!{\bold R}
{\cal H}om(Y; \Omega^{l-j}_{Y/\Delta}(\log D)(-D)
\otimes\pi^*\Omega^{n-l}_{\Delta}, \omega^{\cdot}_Y) 
\cong R^1\pi_!\Omega^j_{Y/\Delta}(\log D)$. 
Therefore we have to prove that 
${\cal E}xtop^{1-n}(\Delta, 
{\bf R}\pi_*(\Omega^{l-j}_{Y/\Delta}(\log D)(-D)
\otimes\pi^*\Omega^{n-l}_{\Delta}), 
\omega^{\cdot}_{\Delta}) = 0$. 

Choose a bounded complex of {\bf FN} free 
${\cal O}_{\Delta}$ modules $L^{\cdot}$ representing  
${\bold R}\pi_*\Omega^{n-j}_{Y/\Delta}(\log D)(-D)$ 
(cf. [R-R]). 
Since $R^i\pi_*\Omega^{n-j}_{Y/\Delta}(\log D)(-D) = 0$ 
for $i \geq l-1$, we have ${\cal H}^i(L^{\cdot}) = 0$ 
for $i \geq l-1$. 
Let $Q := \mathrm{Ker}[L^{l-2} \to L^{l-1}]$. 
Then ${\bold R}{\cal H}omtop(\Delta, 
{\bold R}\pi_*
(\Omega^{l-j}_{Y/\Delta}(\log D)(-D)
\otimes\pi^*\Omega^{n-l}_{\Delta}), 
\omega^{\cdot}_{\Delta})$ is represented by the complex 
${\cal H}omtop( ... \to L^{l-3} 
\to Q \to 0 ..., \Omega^{n-l}_{\Delta}[n-l])$. 
Hence we know that ${\cal E}xtop^{1-n}(\Delta, 
{\bold R}\pi_*(\Omega^{l-j}_{Y/\Delta}(\log D)(-D)
\otimes\pi^*\Omega^{n-l}_{\Delta}), 
\omega^{\cdot}_{\Delta}) = 0$.  

We are now in a position to justify the claim. 
Let $j: Y \setminus D' \to Y$.
The $\omega$ defines an element of 
$(\pi_*j_*j^*\Omega^2_Y(\log D))_0$. 
By a coboundary map $(\pi_*j_*j^*\Omega^2_Y(\log D))_0 
\to (R^1\Gamma_{\pi, D'}\Omega^2_Y(\log D))_0$, 
it defines the obstruction class $ob(\omega) 
\in (R^1\Gamma_{\pi, D'}\Omega^2_Y(\log D))_0$. 
By the observation above, $ob(\omega)$ 
is sent to zero by the natural map 
$(R^1\Gamma_{\pi, D'}\Omega^2_Y(\log D))_0 
\to (R^1\pi_{!}\Omega^2_Y(\log D))_0$. 
Therefore there is a small disc 
$\Delta_{\epsilon} \subset \Delta$ and a subset 
$K \subset \pi^{-1}(\Delta_{\epsilon})(= Y_{\epsilon})$ 
which contains $(\pi\vert_{D'})^{-1}(\Delta_{\epsilon})
( = D'_{\epsilon})$ and which is proper over 
$\Delta_{\epsilon}$, such that $ob(\omega)$ 
is already zero in $H^1_K(Y_{\epsilon}, 
\Omega^2_{Y_{\epsilon}}(\log D'_{\epsilon}))$. 
This is nothing but our claim. \vspace{0.2cm}

(b) $k$: general 

(b-1):  Take a general $l + k$ dimensional complete 
intersection $H := H_1 \cap ... \cap H_{n-k-l}$. 
Let $p \in H \cap f(F)$. $p \in f(F)$ is a smooth point. 
Replace $X$ by a small open neighborhood of $p$. 
Then $H \cap f(F) = \{p\}$. By perturbing $H$, 
we can define a flat holomorphic map 
$g: X \to \Delta^{n-k-l}$ with $g^{-1}(0) = H$. 
We may assume that $g$ has a section passing through 
$p$ and each fiber $g^{-1}(t)$ intersects $f(F)$ 
only in this section. The map $f : Y \to X$ 
gives a simultaneous resolution of $X_t$ 
($t \in \Delta^{n-l-k}$). Since $H$ is general and 
$X$ is sufficiently small, 
$D_t := D \cap Y_t$ are normal crossing divisors of 
$Y_t$ for all $t \in \Delta^{n-l-k}$. 
Let $D'$ be the union of irreducible components of 
$D$ which are mapped in the section. 
Then $D' \to \Delta^{n-k-l}$ is a proper map. 
We put $\pi = g \circ f$. 
We often write $\Delta$ for $\Delta^{n-k-l}$. 
By an induction hypothesis we have an isomorphism 
$H^0(Y \setminus D', \Omega^2_{Y \setminus D'}(\log D)) 
\cong H^0(U, \Omega^2_U)$.  

(b-2): We shall prove the following \vspace{0.15cm}

{\bf Claim}. {\em By replacing $\Delta^{n-l}$ 
by a smaller disc and by restricting everything 
(e.g. $X$, $Y$, $D$, $D'$ ...) over the new disc, 
we have a subset $K \subset Y$ which contains $D'$ 
and which is proper over $\Delta^{n-l}$ 
with the following property:}  

{\em The $\omega$ is mapped to zero by the 
comoposition of the maps}

$$   H^0(U, \Omega^2_U) \cong 
H^0(Y \setminus D', \Omega^2_{Y \setminus D'}(\log D)) 
\to H^0(Y \setminus K, \Omega^2_{Y \setminus K}(\log D)) 
\to H^1_K(Y, \Omega^2_Y(\log D)). $$    \vspace{0.2cm}

If the claim is verified, then 
$\omega\vert_{Y \setminus K}$ extends to a logarithmic 
2-form on $Y$. It is clear that its restriction 
to $U$ is $\omega$.  

The proof of the claim is similar to the claim in (a-2). 
When we apply the relative duality we need the vanishings: 
$R^i\pi_*\Omega^{l+k-2}_{Y_t}(\log D_t)(-D_t) 
= R^i\pi_*\Omega^{l+k-1}_{Y_t}(\log D_t)(-D_t) 
= R^i\pi_*\omega_{Y_t} = 0$ 
for $i \geq l+k-1$ and for $t \in \Delta$. 
\hspace{1.0cm} Q.E.D.  \vspace{0.2cm}

{\bf Lemma 2}. {\em Let $p \in X$ be a Stein open 
neighborhood of a point $p$ of a complex algebraic 
variety. Assume that 
$X$ is a rational singularity of 
$dim X \geq 3$. Let $f : Y \to X$ be a resolution of 
singularities of $X$ such that $E := f^{-1}(p)$ is a 
simple normal crossing divisor. Then 
$H^0(Y, \Omega^i_Y) \to H^0(Y, \Omega^i_Y(\log E))$ 
are isomorphisms for $i = 1, 2$.} \vspace{0.15cm}

{\em Proof}. By the assumption we can take a complete 
algebraic variety $Z$ which contains $X$ as an open 
set. We may assume that $f$ is obtained from a 
resolution $\tilde Z \to Z$. Set $V := \tilde Z 
\setminus E$. Recall that the 
natural exact sequence  

$$ \to H^{j}({\tilde Z}, {\bold C}) \to 
H^{j}(V, {\bold C}) \to 
H^{j+1}_E({\tilde Z}, {\bold C}) \to $$ 
is obtained from the following exact 
sequence of the complexes by taking 
 hypercohomology 

$$ 0 \to \Omega^{\cdot}_{\tilde Z} 
\to \Omega^{\cdot}_{\tilde Z}(\log E) 
\to \Omega^{\cdot}_{\tilde Z}(\log E)/
\Omega^{\cdot}_{\tilde Z} \to 0. $$ 

Introduce the stupid filtrations $F^{\cdot}$
(cf. [De]) on three complexes and take 
${\bold H}^j(Gr^i_{F})$ of 
the sequence of complexes. 
Then we have 

$$ \to H^{j-i}(\Omega^i_{\tilde Z}) \to 
H^{j-i}(\Omega^i_{\tilde Z}(\log E)) \to 
H^{j-i}(\Omega^i_{\tilde Z}(\log E)/\Omega^i
_{\tilde Z}) \to . $$   

We know that this exact sequence coincides 
with the exact sequence 

$$ \to Gr^i_F(H^j(\tilde Z, {\bold C}) 
\to Gr^i_F(H^j(V, {\bold C}) \to 
Gr^i_F(H^{j+1}_E(\tilde Z, {\bold C})) \to 
$$
which comes from the mixed Hodge structures. 
In particular, the map 
$H^{j-i}(\Omega^i_{\tilde Z}(\log E)/\Omega^i
_{\tilde Z}) \to H^{j-i+1}(\Omega^i_{\tilde Z})$ 
is interpreted as the map $Gr^i_F(H^{j+1}_E
({\tilde Z}, {\bold C})) \to Gr^i_F(H^{j+1}
({\tilde Z}, {\bold C}))$.    

We next consider the natural map of 
mixed Hodge structures: 
$H^{j+1}({\tilde Z}, {\bold C}) \to 
H^{j+1}(E, {\bold C})$. We have 
$Gr^i_F(H^{j+1}(E, {\bold C})) 
\cong H^{j-i+1}(\hat\Omega^i_E)$ 
(cf. [Fr]), (see Introduction, 
for the notation $\hat\Omega^i_E$). Note 
that $\hat\Omega^i_E$ is 
isomorphic to the cokernel of the 
injection $\Omega^i_Y(\log E)(-E) \to 
\Omega^i_Y$. Therefore 
the composed map 

$$ H^{j-i}(\Omega^i_{\tilde Z}(\log E)/\Omega^i
_{\tilde Z}) \to H^{j-i+1}(\Omega^i_{\tilde Z}) 
\to H^{j-i+1}
(\Omega^i_Y/\Omega^i_Y(\log E)(-E)) $$ 
is interpreted as the map 

$$ Gr^i_F(H^{j+1}_E
({\tilde Z}, {\bold C})) \to 
Gr^i_F(H^{j+1}(E, {\bold C}). $$ 
 
By the isomorphisms 
$H^{\cdot}(Y, {\bold C}) \cong 
H^{\cdot}(E, {\bold C})$, 
$H^{\cdot}(Y)$ have natural 
mixed Hodge structures. 
We put $U' := Y \setminus E$. 

($i = 2$): We shall first prove that the natural map 
$\beta : H^3_E(Y, {\bold C}) \to H^3(Y, {\bold C})$ 
is an injection.

By a local cohomology exact sequence 
it suffices to show that $\alpha : H^2(Y, {\bold C}) 
\to H^2(U', {\bold C})$ is a surjection. Because $X$ 
has only rational singularity, 
$H^1(Y, {\cal O}^{*}_Y)\otimes{\bold C} 
\cong H^2(Y, {\bold C})$.

 On the other hand, one can prove that 
$H^1(U', {\cal O}^{*}_{U'})\otimes{\bold C} 
\to H^2(U', {\bold C})$ is a surjection; 
in fact, since $H^{\cdot}(U', {\bold C}) 
\cong {\bold H}^{\cdot}
(Y, {\bold R}j_*\Omega^{\cdot}_{U'}) 
\cong {\bold H}^{\cdot}(Y, j_*\Omega^{\cdot}_{U'}) 
\cong {\bold H}^{\cdot}(Y, \Omega^{\cdot}_Y(\log E))$ 
where $j : U' \to Y$ is the immersion [De], 
and since $H^{\cdot}(U', {\bold C}) 
\cong H^{\cdot}(U', \Omega^{\cdot}_{U'})$, 
we have a commutative diagram of Hodge spectral sequences

\begin{equation}
\begin{CD}
H^q(Y, \Omega^p_Y(\log E))  @=> H^{p+q}(U', {\bold C}) \\
@VVV @VVV \\
H^q(U', \Omega^p_{U'}) @=>  H^{p+q}(U', {\bold C})  
\end{CD}
\end{equation}

Let $F^{\cdot}_1$ (resp. $F^{\cdot}_2$) be the filtration 
of $H^{p+q}(U', {\bold C})$ by the first (resp. second) 
spectral sequence. In particular, when $p + q = 2$, 
we have a surjection $Gr^0_{F_1}H^2(U', {\bold C}) 
\to Gr^0_{F_2}H^2(U', {\bold C})$. Because $X$ has 
only rational singularity, $H^2(Y, {\cal O}_Y) = 0$, 
hence $Gr^0_{F_1}H^2(U', {\bold C}) = 0$. 
Therefore $Gr^0_{F_2}H^2(U', {\bold C}) = 0$, 
which implies that the natural map 
$H^2(U', {\bold C}) \to H^2(U', {\cal O}_{U'})$ 
is the zero map. It is now clear that 
$H^2(U', {\bold Z}) \to H^2(U', {\cal O}_{U'})$ 
is also the zero map, and hence 
$H^1(U', {\cal O}^{*}_{U'}) \to H^2(U', {\bold Z})$ 
is a surjection. Therefore one has a surjection 
$H^1(U', {\cal O}^{*}_{U'})\otimes{\bold C} 
\to H^2(U', {\bold C})$.     

It is enough to prove that $H^1(Y, {\cal O}^*_Y) 
\to H^1(U', {\cal O}^*_{U'})$ 
is surjective in order to prove that 
$\alpha$ is surjective. Put $f^0 := f\vert_{U'}: U' 
\to U$, where $U := X \setminus \{p\}$. 
Let $D_i$ be irreducible exceptional divisors of 
$f$ where $f(D_i) \neq \{p\}$. Let 
$L \in \mathrm{Pic}(U')$. 
Since $\dim X \geq 3$ and $U \subset X$ is 
1-concave at $p$, there exists a reflexive coherent sheaf 
$F$ on $X$ of rank 1 such that 
$F\vert_U \cong (f^0_*L)^{**}$ by [S, Theorem 5], 
where $**$ means the double dual. 
By taking the double dual of both sides 
of the natural map 
$(f^0)^*f^0_*L \to L$ we get an injection 
$((f^0)^*f^0_*L)^{**} \to L$, hence 
$L \cong ((f^0)^*f^0_*L)^{**}
\otimes{\cal O}_{U'}(\Sigma a_iD_i)$ 
with some $a_i \geq 0$. 
On the other hand, we have an injection 
$((f^0)^*f^0_*L)^{**} \to (f^*f_*F)^{**}\vert_{U'}$ 
and $M: = (f^*f_*F)^{**} \in \mathrm{Pic}(Y)$. 
We have $((f^0)^*f^0_*L)^{**} 
\cong M\vert_{U'} \otimes {\cal O}_{U'}(\Sigma (-b_i)D_i)$ 
with some $b_i > 0$. Therefore 
$L \cong M\otimes 
{\cal O}_Y(\Sigma (a_i - b_i)D_i)\vert_{U'}$.  
\vspace{0.15cm}
 
We are now in a position to prove the 
lemma for $(i = 2)$. Let us consider 
the exact sequence

$$ 0 \to \Omega^2_Y/\Omega^2_Y(\log E)(-E) 
\to \Omega^2_Y(\log E)/\Omega^2(\log E)(-E) 
\to \Omega^2_Y(\log E)/\Omega^2_Y \to 0. $$

From this sequence we have a map $\delta : 
H^0(E, \Omega^2_Y(\log E)/\Omega^2_Y) 
\to H^1(E, \Omega^2_Y/\Omega^2_Y(\log E)(-E))$. 
The map $\beta$ is a morphism of mixed Hodge structures 
and $\delta$ can be interpreted as the map 
$Gr^2_F(H^3_E(Y, {\bold C}) \to Gr^2_F(H^3(Y, {\bold C})$. 
We already proved that $\beta$ is an injection. 
Hence $\delta$ is also an injection by the 
strict compatibility of the filtrations $F$. 
Note that $\delta$ is factorized as 
$H^0(E, \Omega^2_Y(\log E)/\Omega^2_Y)
 \stackrel{\gamma}\to H^1(Y, \Omega^2_Y) 
\to H^1(E, \Omega^2_Y/\Omega^2_Y(\log E)(-E))$ 
where $\gamma$ is the last map in the 
following exact sequence

$$  H^0(Y, \Omega^2_Y) \stackrel{\tau}\to 
H^0(Y, \Omega^2(\log E)) \to 
H^0(E, \Omega^2_Y(\log E)/\Omega^2_Y) 
\to H^1(Y, \Omega^2_Y). $$

Since $\delta$ is injective, $\gamma$ is also injective.  
Hence $\tau$ is surjective by the exact sequence.

($i = 1$): We shall first prove that the natural map 
$\beta : H^2_E(Y, {\bold C}) \to H^2(Y, {\bold C})$ 
is an injection. By a local cohomology exact sequence 
it suffices to show that $\alpha : H^1(Y, {\bold C}) 
\to H^1(U', {\bold C})$ is a surjection. Since $X$ has
 rational singularities, the sequence 

$$  H^0(Y, {\cal O}_Y) \to H^0(Y, {\cal O}^*_Y) 
\to H^1(Y, {\bold Z}) \to 0 $$
is exact. By the same argument as ($i = 2$), 
$H^1(U', {\bold Z}) \to H^1(U', {\cal O}_{U'})$ 
is the zero map because $X$ has rational singularities. 
Hence the sequence 

$$  H^0(U', {\cal O}_{U'}) \to H^0(U', {\cal O}^*_{U'}) 
\to H^1(U', {\bold Z}) \to 0  $$

is also exact. The restriction map 
$H^0(Y, {\cal O}^*_Y) \to H^0(U', {\cal O}^*_{U'})$ 
is an isomorphism because it is factorized as 
$H^0(Y, {\cal O}^*_Y) \cong H^0(X, {\cal O}^*_X) 
\cong H^0(X\setminus \{p\}, 
{\cal O}^*_{X \setminus \{p\}}) 
\cong H^0(U', {\cal O}^*_{U'})$. 
Similarly the restriction map 
$H^0(Y, {\cal O}_Y) \to H^0(U', {\cal O}_{U'})$ 
is also an isomorphism. 
Hence the restriction $H^1(Y, {\bold Z}) 
\to H^1(U', {\bold Z})$ is an isomorphism 
by the exact sequences above, and 
$\alpha$ is an isomorphism.  
 
Let us consider the exact sequence 

$$ 0 \to \Omega^1_Y/\Omega^1_Y(\log E)(-E) 
\to \Omega^1_Y(\log E)/\Omega^1(\log E)(-E) 
\to \Omega^1_Y(\log E)/\Omega^1_Y \to 0. $$

From this sequence we have a map 
$\delta : H^0(E, \Omega^1_Y(\log E)/\Omega^1_Y) 
\to H^1(E, \Omega^1_Y/\Omega^1_Y(\log E)(-E))$. 
The map $\beta$ is a morphism of mixed Hodge structures 
and $\delta$ can be interpreted as the map 
$Gr^1_F(H^2_E(Y, {\bold C}) 
\to Gr^1_F(H^2(Y, {\bold C})$. 
We already proved that $\beta$ is an injection. 
Hence $\delta$ is also an injection by the 
strict compatibility of the filtrations $F$. 
Note that $\delta$ is factorized as 
$H^0(E, \Omega^1_Y(\log E)/\Omega^1_Y) 
\stackrel{\gamma}\to H^1(Y, \Omega^1_Y) 
\to H^1(E, \Omega^1_Y/\Omega^1_Y(\log E)(-E))$ 
where $\gamma$ is the last map in 
the following exact sequence

$$  H^0(Y, \Omega^1_Y) \stackrel{\tau}\to 
H^0(Y, \Omega^1(\log E)) 
\to H^0(E, \Omega^1_Y(\log E)/\Omega^1_Y) 
\to H^1(Y, \Omega^1_Y). $$

Since $\delta$ is injective, $\gamma$ is also injective.  
Hence $\tau$ is surjective by the exact sequence.  
Q.E.D.  \vspace{0.15cm}

{\bf Remark}. In the proof of Lemma 2 the map 
$H^0(E, \Omega^i_Y/\Omega^1(\log E)(-E)) 
\to H^0(E, \Omega^i_Y(\log E)/\Omega^i_Y(\log E)(-E))$ 
is surjective for $i = 1, 2$  because we have proved that 
$\delta$ is injective. \vspace{0.15cm}

{\bf Proposition 3}. {\em Let $X$ be a Stein open subset
 of a complex algebraic variety. 
Assume that $X$ has only rational Gorenstein 
singularities. 
Let $\Sigma$ be the singular locus of $X$ and 
let $f : Y \to X$ be a resolution of singularities 
such that $D := f^{-1}(\Sigma)$ is a simple normal 
crossing divisor and such that 
$f\vert_{Y \setminus D} : Y \setminus D 
\cong X \setminus \Sigma$. 
Then $f_*\Omega^2_Y \cong f_*\Omega^2_Y(\log D)$.}  
\vspace{0.2cm}

{\em Proof}. Since $(X, p) \cong (R.D.P.)\times 
({\bold C}^{n-2}, 0)$ 
for all $p \in X$ outside certain codimension 3 
(in $X$) locus 
$\Sigma_0 \subset \Sigma$ ([Re]), it is clear that 
$f_*\Omega^2_Y 
\cong f_*\Omega^2_Y(\log D)$ at such points $p$.   

Let $\omega \in H^0(Y, \Omega^2_Y(\log D))$. 
Let $F$ be an irreducible 
component of $D$ with $f(F) \subset \Sigma_0$. 
Put $k := \dim \Sigma_0 - \dim f(F)$. 
We shall prove that $\omega$ is regular along $F$ by the 
induction on $k$.  

(a) $k = 0$: 

(a-1): Put $l := \mathrm{codim}(\Sigma_0 \subset X)$. 
Note that $l \geq 3$. 
Take a general $l$ dimensional complete intersection 
$H := H_1 \cap H_2 \cap ... \cap H_{n-l}$. 
Let $p \in H \cap f(F)$. $p \in f(F)$ is a smooth point. 
Replace $X$ by a small open neighborhood of $p$. 
Then $H \cap f(F) = \{p\}$. Moreover, 
$\tilde H := f^{-1}(H)$ is a resolution of 
singularities of $H$. Since $X$ has canonical 
singularities, 
$H$ has also canonical singularities. 
$H$ has a unique dissident point $p$ and 
other singular points are locally isomorphic to 
$(R.D.P.) \times ({\bold C}^{l-2}, 0)$. 
By perturbing $H$ we can define a flat 
holomorphic map $g : X \to \Delta^{n-l}$ with 
$g^{-1}(0) = H$. 
We may assume that $g$ has a section passing through 
$p$ and each fiber $X_t := g^{-1}(t)$ intersects $f(F)$ 
only in this section. 
Denote by $p_t \in X_t$ this intersection point. 
By definition $p_0 = p$. The map $f : Y \to X$ 
gives a simultaneous resolution of $X_t$ for $t \in 
\Delta^{n-l}$. 
Let $D'$ be the union of irreducible components of $D$ 
which are mapped in this section. 
Since $H$ is general and $X$ is sufficiently small, 
every irreducible component of $D'$ is mapped 
onto the section. $D' \to \Delta^{n-l}$ 
is a proper map and every fiber $D'_t$ is a 
normal crossing variety.  
Note that $f_t^{-1}(p_t) = D'_t$.  
We put $\pi = g\circ f$. 
We often write $\Delta$ for $\Delta^{n-l}$. 
There are filtrations 
$(\pi\vert_{D'})^*\Omega^2_{\Delta} \subset 
{\cal F} \subset \hat\Omega^2_{D'}$ and 
$\pi^*\Omega^2_{\Delta} \subset {\cal G} \subset 
\Omega^2_Y(\log D')$ which yield the following 
exact sequences 

$$ 0 \to {\cal F} \to \hat\Omega^2_{D'} 
\to \hat\Omega^2_{D'/\Delta} 
\to 0, $$

$$ 0 \to (\pi\vert_{D'})^*\Omega^2_{\Delta} 
\to {\cal F} \to 
(\pi\vert_{D'})^*\Omega^1_{\Delta}\otimes 
\hat\Omega^1_{D'/\Delta} \to 0, 
$$ 

$$ 0 \to {\cal G} \to \Omega^2_Y(\log D') 
\to \Omega^2_{Y/\Delta}(\log 
D') \to 0, $$

$$ 0 \to \pi^*\Omega^2_{\Delta} \to 
{\cal G} \to 
\pi^*\Omega^1_{\Delta}\otimes 
\Omega^1_{Y/\Delta}(\log D') \to 0. $$
  
(a-2): Let us consider the exact sequence 

$$ 0 \to \Omega^2_Y/\Omega^2_Y(\log D')(-D') 
\to \Omega^2_Y(\log D')/\Omega^2(\log D')(-D') 
\to \Omega^2_Y(\log D')/\Omega^2_Y \to 0. $$

We shall prove the following. \vspace{0.15cm}

{\bf Claim}. {\em The map $H^0(D', 
\Omega^2_Y/\Omega^2(\log D')(-D')) 
\to H^0(D', 
\Omega^2_Y(\log D')/\Omega^2_Y(\log D')(-D'))$ 
is surjective.}
 \vspace{0.15cm}

{\em Proof}. Note that 
$\Omega^2_Y/\Omega^2_Y(\log D')(-D') 
\cong \hat\Omega^2_{D'}$. By the exact sequences above 
we have 
two commutative diagrams with exact columns 

\begin{equation}
\begin{CD}
0  @. 0 \\
@VVV  @VVV \\
H^0(D', {\cal F}) @>>>  H^0(D', {\cal G}
\otimes{\cal O}_{D'}) \\
@VVV  @VVV \\
H^0(D', \hat\Omega^2_{D'}) @>>> 
H^0(D', \Omega^2_Y(\log D')\otimes{\cal O}_{D'}) \\
@VVV  @VVV \\
H^0(D', \hat\Omega^2_{D'/\Delta}) @>{\mu_2}>> 
H^0(D', \Omega^2_{Y/\Delta}(\log D')\otimes{\cal O}_{D'}) 
\end{CD}
\end{equation}

\begin{equation}
\begin{CD}
0 @. 0 \\ 
@VVV  @VVV \\
H^0(\Delta, \Omega^2_{\Delta}) @>>>  
H^0(\Delta, \Omega^2_{\Delta}) \\
@VVV  @VVV \\
H^0(D', {\cal F}) @>>> H^0(D, {\cal G}
\otimes{\cal O}_{D'}) \\
@VVV  @VVV \\
H^0(D', \hat\Omega^1_{D'/\Delta})^{\oplus n-l}  
@>{(\mu_1)^{\oplus n-l}}>>  
H^0(D', \Omega^1_{Y/\Delta}(\log D')
\otimes{\cal O}_{D'})^{\oplus n-l}  
\end{CD}
\end{equation}

We shall show that $H^0(D', \hat\Omega^i_{D'/\Delta}) 
= H^0(D', 
\Omega^i_{Y/\Delta}(\log D')\otimes{\cal O}_{D'}) 
= 0$ for $i = 1, 2$. After these are proved our claim 
easily follows from the commutative diagrams above. 

First note that $H^0(D'_t, \hat\Omega^i_{D'_t}) = 0$ 
for $t \in \Delta$ and for $i = 1, 2$. 
In fact, if $h^0(D'_t, \hat\Omega^i_{D'_t}) > 0$, 
then $h^i(D'_t, {\cal O}_{D'_t}) > 0$ 
by the mixed Hodge structures on 
$H^i(D'_t, {\bold C})$. 
On the other hand, $X_t$ has 
only canonical singularities, 
hence has only rational singularities. 
Therefore $H^i(Y_t, {\cal O}_{Y_t}) = 0$. 
Then one can check that 
$H^i(D'_t, {\cal O}_{D'_t}) = 0$. 
This is a contradiction. For details of this argument, 
see [Na 1, Claim 1, (i) in the proof of Prop. (1.1)]. 
As a consequence we know that 
$H^0(D', \hat\Omega^i_{D'/\Delta}) = 0$ 
for $i = 1, 2$.  

We next note that 
$\mu_i(t) : H^0(D'_t, \hat\Omega^i_{D'_t}) 
\to H^0(D'_t, \Omega^i_{Y_t}(\log D'_t)
\otimes{\cal O}_{D'_t})$ are surjective for $i = 1,2$ 
by Remark below Lemma 2. 
Since $H^0(D'_t, \hat\Omega^i_{D'_t}) = 0$, 
this implies that 
$H^0(D'_t, \Omega^i_{Y_t}(\log D'_t)
\otimes{\cal O}_{D'_t}) = 0$, 
hence $H^0(D', \Omega^i_{Y/\Delta}(\log D')
\otimes{\cal O}_{D'}) = 0$ for $i = 1, 2$. 
Q.E.D. \vspace{0.15cm}

(a-3): We shall continue the proof 
in the case $k = 0$. By taking cohomology of 
the exact sequence

$$ 0 \to \Omega^2_Y/\Omega^2_Y(\log D')(-D') 
\to \Omega^2_Y(\log D')/\Omega^2(\log D')(-D') 
\to \Omega^2_Y(\log D')/\Omega^2_Y \to 0. $$
and by applying Claim in (a-2) we see that 
$\delta : H^0(D', \Omega^2_Y(\log D')/\Omega^2_Y) 
\to H^1(D', \Omega^2_Y/\Omega^2_Y(\log D')(-D'))$ 
is an injection. The map $\delta$ is factorized as 
$H^0(D', \Omega^2_Y(\log D')/\Omega^2_Y) 
\stackrel{\gamma}\to H^1(Y, \Omega^2_Y) 
\to H^1(D', \Omega^2_Y/\Omega^2_Y(\log D')(-D'))$ 
where $\gamma$ is the last map 
in the following exact sequence

$$  H^0(Y, \Omega^2_Y) \stackrel{\tau}\to 
H^0(Y, \Omega^2(\log D')) \to 
H^0(D', \Omega^2_Y(\log D')/\Omega^2_Y) 
\to H^1(Y, \Omega^2_Y). $$

Since $\delta$ is injective, $\gamma$ is also injective.  
Hence $\tau$ is surjective by the exact sequence.  
Q.E.D.  \vspace{0.15cm}
   
(b) $k$ : general

(b-1): Take a general $l + k$ dimensional complete 
intersection $H := H_1 \cap H_2 \cap ... \cap H_{n-l-k}$. 
Let $p \in H \cap f(F)$. $p \in f(F)$ is a smooth point. 
Replace $X$ by a small open neighborhood of $p$. 
Then $H \cap f(F) = \{p\}$. Moreover, 
$\tilde H := f^{-1}(H)$ is a resolution of 
singularities of $H$. 
Since $X$ has canonical singularities, $H$ has 
also canonical singularities. 
By perturbing $H$ we can define a flat 
holomorphic map $g : X \to \Delta^{n-l-k}$ with 
$g^{-1}(0) = H$. We may assume that $g$ has a section 
passing through $p$ and each fiber $X_t := g^{-1}(t)$ 
intersects $f(F)$ only in this section. 
Denote by $p_t \in X_t$ this intersection point. 
By definition $p_0 = p$. The map $f : Y \to X$ 
gives a simultaneous resolution of $X_t$ for $t \in 
\Delta^{n-l-k}$. 
Let $D'$ be the union of irreducible components of $D$ 
which are mapped in this section. 
Since $H$ is general and $X$ is sufficiently small, 
every irreducible component of $D'$ 
is mapped onto the section. 
$D' \to \Delta^{n-l}$ is a proper map and 
every fiber $D'_t$ is a normal crossing 
variety.  Note that $f_t^{-1}(p_t) = D'_t$.  
We put $\pi = g\circ f$. 
We often write $\Delta$ for $\Delta^{n-l-k}$. 
There are filtrations 
$(\pi\vert_{D'})^*\Omega^2_{\Delta} \subset 
{\cal F} \subset \hat\Omega^2_{D'}$ and 
$\pi^*\Omega^2_{\Delta} \subset \cal G \subset 
\Omega^2_Y(\log D')$ which yield the same exact 
sequences as (a-1).

By an induction hypothesis we have an isomorphism 
$H^0(Y, \Omega^2_Y(\log D') \cong 
H^0(Y, \Omega^2_Y(\log D)$. 
Therefore we have to prove that 
$H^0(Y, \Omega^2_Y) \to H^0(Y, \Omega^2_Y(\log D')$ 
is a surjection (hence an isomorphism). \vspace{0.15cm} 

(b-2): Let us consider the exact sequence 

$$ 0 \to \Omega^2_Y/\Omega^2_Y(\log D')(-D') 
\to \Omega^2_Y(\log D')/\Omega^2(\log D')(-D') 
\to \Omega^2_Y(\log D')/\Omega^2_Y \to 0. $$

We shall prove the following. \vspace{0.15cm}

{\bf Claim}. {\em The map $H^0(D', 
\Omega^2_Y/\Omega^2(\log D')(-D')) 
\to H^0(D', 
\Omega^2_Y(\log D')/\Omega^2_Y(\log D')(-D'))$ 
is surjective.}
 \vspace{0.15cm}

The proof is similar to Claim in (a-2). \vspace{0.15cm}

(b-3):  Once Claim is proved, then 
$H^0(Y, \Omega^2_Y) \to H^0(Y, \Omega^2_Y(\log D')$ 
is an isomorphism by the same argument as (a-3).  Q.E.D.  
\vspace{0.15cm}

By combining Propositions 1 and 3 
we have the following. \vspace{0.15cm}

{\bf Theorem 4.}. {\em  Let $X$ be a Stein open subset
 of a complex algebraic variety. 
Assume that $X$ has only rational 
Gorenstein singularities. 
Let $\Sigma$ be the singular locus of $X$ and 
let $f : Y \to X$ be a resolution of 
singularities such that 
$f\vert_{Y\setminus f^{-1}(\Sigma)} : 
Y \setminus f^{-1}(\Sigma) \cong X \setminus \Sigma$. 
Then $f_*\Omega^2_Y \cong i_*\Omega^2_U$ 
where $U := X \setminus \Sigma$ and $i : U \to X$ 
is a natural injection.} \vspace{0.15cm}

\begin{center}
{\bf \S 2. Symplectic Varieties}
\end{center}
\vspace{0.15cm}

We shall begin this section with the stability of 
Kaehlerity under small deformation. \vspace{0.15cm}

{\bf Proposition 5}. {\em Let $Z$ be a compact normal 
Kaehler space with rational singularities. 
Then any small (flat) deformation $Z_t$ of $Z$ is also a 
Kaehler space.}  \vspace{0.15cm}

{\em Proof}. By a theorem of Bingener [B] 
we only have to prove that the map 
$H^2(Z, {\bold R}) \to H^2(Z, {\cal O}_Z)$ 
induced by a sheaf homomorphism 
${\bold R}_Z \to {\cal O}_Z$ is surjective. 
Let $f : \tilde Z \to Z$ be a 
resolution of singularities. 
Since $Z$ has only rational singularities, 
we know that $R^1f_*{\bold R}_{\tilde Z} = 0$. 
Consider the commutative diagram with exact rows 

\begin{equation}
\begin{CD}
0 @>>> H^2(Z, {\bold R}) @>>>  H^2(\tilde Z, {\bold R}) 
@>>> H^0(Z, R^2f_*{\bold R}) \\ 
@. @VVV  @VVV @VVV \\
0 @>>> H^2(Z, {\cal O}_Z) @>>> 
H^2(\tilde Z, {\cal O}_{\tilde Z}) @>>> 0 
\end{CD}
\end{equation}

By Hodge theory the middle vertical map is surjective. 
The surjectivity of the left vertical map follows 
from the following claim. \vspace{0.15cm}

{\bf Claim}. $\mathrm{im}[H^2(\tilde Z, {\bold R}) 
\stackrel{\phi}\to H^0(Z, R^2f_*{\bold R})] 
= \mathrm{im}[H^2(\tilde Z, {\bold R})\cap H^{1,1} 
\to H^0(Z, R^2f_*{\bold R})]$.   \vspace{0.15cm}

{\em Proof}. It suffices to show that, if an element 
$\alpha \in H^2(\tilde Z, {\bold R})$ has 
the Hodge decomposition $\alpha = 
\alpha^{(2,0)} + \alpha^{(0,2)}$, then 
$\phi(\alpha) = 0$.  
Since $\bar\alpha^{(2,0)} = \alpha^{(0,2)}$, 
we only have to prove that 
$\phi_{\bold C}(\alpha^{(2,0)}) = 0$. 
We shall show that, for every point $z \in Z$, 
$\phi_{\bold C}(\alpha^{(2,0)})_z = 0$ 
in $(R^2f_*{\bold C})_z$. 
Let $\nu : W \to \tilde Z$ be a projective 
bimeromorphic map such that $W$ is smooth and 
$D := (f\circ\nu)^{-1}(x)$ is a simple normal 
crossing divisor of $W$. 
Put $h := f\circ \nu$. 
Since $R^1\nu_*{\bold C} = 0$, 
$R^2f_*{\bold C}$ injects to $R^2h_*{\bold C}$. 
Therefore we have to check that $\alpha^{(2,0)}$ 
is sent to zero by the composition of the maps 
$H^2(\tilde Z, {\bold C}) \to H^2(W, {\bold C}) 
\to (R^2h_*{\bold C})_z (= H^2(D, {\bold C}))$. 
Since the map $H^2(\tilde Z, {\bold C}) 
\to H^2(D, {\bold C})$ 
is a morphism of mixed Hodge structures, 
it preserves the Hodge filtration $F$. 
In particular, it induces 
$Gr^2_F(H^2(\tilde Z, {\bold C})) 
\to Gr^2_F(H^2(D, {\bold C}))$. 
Since $\alpha^{(2,0)} \in 
Gr^2_F(H^2(\tilde Z, {\bold C}))$, 
$\phi_{\bold C}(\alpha^{(2,0)}) 
\in Gr^2_F(H^2(D, {\bold C}))$. 
On the other hand, $Gr^2_F(H^2(D, {\bold C})) = 0$ 
because $Z$ has rational singularities; 
if $Gr^2_F(H^2(D, {\bold C})) \neq 0$, 
then $Gr^0_F(H^2(D, {\bold C})) 
= H^2(D, {\cal O}_D) \neq 0$ which  contradicts 
the fact that $Z$ has rational singularities. 
For details of the argument, see 
[Na, Claim 1, (i) in the proof of Prop. (1.1)].   
Q.E.D.  \vspace{0.15cm}

By Theorem 4 we have another definition of a 
symplectic singularity. \vspace{0.15cm}

{\bf Theorem 6}. {\em Let $X$ be a Stein open subset 
of a complex algebraic variety. 
Then $X$ is a symplectic singularity 
(for the definition, see Introcuction) 
if and only if $X$ has rational Gorenstein 
singularities and the regular locus  
$U$ of $X$ admits a everywhere non-degenerate 
holomorphic closed 2-form.}  \vspace{0.15cm}

{\em Proof}. By [Be 1] if $X$ is symplectic, 
then $X$ has rational Gorenstein singularities. 
Hence, "only if"part holds. 
On the other hand, by Theorem 4, "if" part holds. 
Q.E.D.    \vspace{0.15cm}

{\bf Example 6'}(cf. [O, 1.5]): 
Let $c$ be an even positive integer with
$c \geq 4$ and 
let $E$ be a ${\bold C}$ vector space with 
$\dim E = c$, equipped with a non-degenerate 
(alternative) 2-form $\omega : E \times E \to {\bold C}$. 
Let $W$ be a 3-dimensional ${\bold C}$ 
vector space with a non-degenerate symmetric form 
$\kappa : W \times W \to {\bold C}$. 
Let $SO(W)$ be the special orthogonal subgroup of 
$GL(W)$ with respect to $\kappa$. Put 

$$   \mathrm{Hom}^{\omega}(W, E) 
:= \{ \phi \in \mathrm{Hom}(W, E); 
\phi^*\omega = 0 \}. $$

Let $X := \mathrm{Hom}^{\omega}(W, E)//SO(W)$. 
Identifying $\mathrm{Hom}(W, E) = W^* 
\otimes E$ with $W \otimes E$ by $\kappa$, 
we can define a 2-form $\tilde\omega$ on 
$\mathrm{Hom}(W, E)$ by $\tilde\omega(\alpha 
\otimes \beta, \alpha'\otimes \beta') 
:= \kappa(\alpha, \alpha')\omega(\beta, \beta')$ 
for $\alpha, \alpha' \in W$, $\beta, \beta' \in E$. 
Let $\mathrm{Hom}^{\omega}(W, E)^s$ be the 
open subset of $\mathrm{Hom}^{\omega}(W, E)$
 which consists of points 
with trivial isotropy group and with closed 
orbits. 
Then $U:= \mathrm{Hom}^{\omega}(W, E)^s/SO(W)$ 
becomes the regular part of $X$. 
Put $\Sigma := \mathrm{Sing}(X)$. 
By calculations $\dim X = 3c - 6$ and 
$\dim \Sigma = c$. The 2-form 
$\tilde\omega\vert_{\mathrm{Hom}^{\omega}(W, E)^s}$
 descends to a symplectic 2-form $\omega_U$ on $U$. 
By a theorem of Boutot [Bo] we see that $X$ has 
rational singularities because 
$\mathrm{Hom}^{\omega}(W, E)$ 
has only rational singualrities. 
$\wedge^{(3/2)c - 3}\omega_U$ 
gives a trivialization of the dualizing sheaf of 
$X$; hence $X$ has rational 
Gorenstein singularities. By Theorem 6 $X$ 
has symplectic singualrities. \vspace{0.2cm}   

{\bf Theorem 7}. {\em Let $(Z, \omega)$ be a 
projective symplectic variety.  
Let $g : {\cal Z} \to \Delta$ be a projective
 flat morphism from ${\cal Z}$ to a 
1-dimensional unit disc 
$\Delta$ with $g^{-1}(0) = Z$. 
Then $\omega$ extends sideways 
in the flat family so that 
it gives a symplectic 2-form $\omega_t$
 on each fiber $Z_t$ for 
$t \in \Delta_{\epsilon}$ 
with a sufficiently small $\epsilon$.}  
\vspace{0.15cm}

{\em Proof}. We shall shrink $\Delta$ suitably 
in each step of arguments, 
but use the same notation $\Delta$ after shrinking. 
Put $\dim Z = 2l$.

Let $\nu_0: \tilde Z \to Z$ be a resolution. 
Let $\omega \in 
H^0(Z, {\nu_0}_*\Omega^2_{\tilde Z})$ 
be a symplectic 2-form. We take the conjugate 
$\bar{\omega} \in H^2(\tilde Z, {\cal O}_{\tilde Z})$ 
by the Hodge decomposition 
$H^2(\tilde Z, {\bold C}) = 
H^0(\tilde Z, \Omega^2_{\tilde Z})\oplus 
H^1(\tilde Z, \Omega^1_{\tilde Z})\oplus 
H^2(\tilde Z, {\cal O}_{\tilde Z})$. 
Since $Z$ has only rational singularities, 
$H^2(\tilde Z, {\cal O}_{\tilde Z}) \cong 
H^2(Z, {\cal O}_Z)$; by this isomorphism 
we regard $\bar{\omega}$ as an element of 
$H^2(Z, {\cal O}_Z)$. Since $Z_t$ are projective 
varieties with rational singularities, 
the natural maps $H^i(Z_t, {\bold C}) 
\to H^i(Z_t, {\cal O}_{Z_t})$ are surjective 
for all $i$ by [Ko, Theorem 12.3]; hence by [D-J] 
$R^ig_*{\cal O}_{\cal Z}$ are locally free sheaves 
which are compatible with base change. 
Therefore $\bar{\omega}$ extends sideways and 
defines non-zero $\bar{\omega_t} 
\in H^2(Z_t, {\cal O}_{Z_t})$ for each $t$. 
Since $\wedge^l \bar{\omega} \in H^{2l}(Z, {\cal O}_Z) 
= {\bold C}$ is not zero, we also have 
$\wedge^{l}\bar{\omega_t} \neq 0$ in 
$H^{2l}(Z_t, {\cal O}_{Z_t}) = {\bold C}$. 
Take a resolution $\tilde{Z_t} \to Z_t$, 
and identify $H^2(Z_t, {\cal O}_{Z_t})$ 
with $H^2(\tilde{Z_t}, {\cal O}_{\tilde Z_t})$. 
By Hodge decomposition of 
$H^2(\tilde Z_t, {\bold C})$, 
we take the conjugate 
$\omega_t \in H^0(\tilde Z_t, \Omega^2_{\tilde Z_t})$ 
of $\bar{\omega_t}$. Note that 
$0 \neq \wedge^{l}\omega_t \in 
H^0(\tilde Z_t, \omega_{\tilde Z_t}) = {\bold C}$. 
This implies that $\omega_t$ is everywhere 
non-degenerate at regular locus of $Z_t$ 
because $\omega_{Z_t}$ is trivial, 
and we know that $Z_t$ is a symplectic variety. 

However, by the fiberwise argument above, 
it is not clear whether $\omega$ holomorphically 
extends sideways. 
We shall prove that $\omega$ actually extends sideways 
by using Theorem 4. 
Let ${\cal U} := \{z \in {\cal Z}$; 
$g$ is smooth at $z$ $\}$. 
Denote by $i$ the natural inclusion of ${\cal U}$ to 
${\cal Z}$. Put $F := \Omega^2_{{\cal U}/\Delta}$, 
$F_0 := \Omega^2_U$ and $U := {\cal U} \cap Z$. 
$i_*F$ is a coherent torsion free sheaf on 
${\cal Z}$, and hence is flat over $\Delta$. 
By the exact sequence 

$$  0 \to i_*F \stackrel{t}\to i_*F \to i_*F_0   $$ 

we know that $i_*F\otimes_{{\cal O}_{\cal Z}}{\cal O}_Z 
\to i_*F_0$ is an injection, 
hence $h^0(i_*F\otimes_{{\cal O}_{\cal Z}}{\cal O}_Z) 
\leq h^0(i_*F_0)$. 

On the other hand, for general $t$, 
$h^0(i_*F\otimes_{{\cal O}_{\cal Z}}{\cal O}_{Z_t}) 
= h^0(i_*F_t)$. This is proved in the following way. 
Since $t \in \Delta$ is general, 
we may assume that $g: {\cal Z} \to \Delta$ has a 
simultaneous 
resolution $\alpha: \tilde{\cal Z} \to {\cal Z}$ 
if we replace $\Delta$ by a suitable 
open neighborhood of 
$t$. Put $f = g\circ\alpha$. 
We have a commutative diagram:

\begin{equation}
\begin{CD}
\alpha_*\Omega^2_{\tilde Z/\Delta}
\otimes_{{\cal O}_{\cal Z}}{\cal O}_{Z_t} @>>> 
i_*F\otimes_{{\cal O}_{\cal Z}}{\cal O}_{Z_t} \\
@VVV  @VVV \\
{\alpha_t}_*\Omega^2_{\tilde Z_t} @>>> i_*F_t 
\end{CD}
\end{equation}

The horizontal map at the bottom is an isomorphism 
by Theorem 4. We shall prove that 
$H^0(\alpha_*\Omega^2_{\tilde Z/\Delta}
\otimes_{{\cal O}_{\cal Z}}{\cal O}_{Z_t}) 
\to H^0({\alpha_t}_*\Omega^2_{\tilde Z_t})$ 
is surjective; if so, then 
$H^0(i_*F\otimes_{{\cal O}_{\cal Z}}{\cal O}_{Z_t}) 
\to H^0(i_*F_t)$ is also a surjection by the diagram, 
and hence is an isomorphism. 
Now apply base change theorem to 
$(\alpha_*\Omega^2_{\tilde{\cal Z}/\Delta}, g)$ 
and $(\Omega^2_{\tilde{\cal Z}/\Delta}, f)$. 
Then we have a commutative diagram 

\begin{equation}
\begin{CD}
g_*(\alpha_*\Omega^2_{\tilde{\cal Z}/\Delta})
\otimes k(t) 
@>>> H^0(Z_t, \alpha_*\Omega^2_{\tilde{\cal Z}/\Delta}
\otimes{\cal O}_{Z_t}) \\
@VVV @VVV \\
f_*\Omega^2_{\tilde{\cal Z}/\Delta}\otimes k(t) 
@>>> H^0(Z_t, {\alpha_t}_*\Omega^2_{\tilde Z_t}) 
\end{CD}
\end{equation}

The vertical map on the left hand side is 
clearly an isomorphism. 
The horizontal maps are both isomorphisms 
if $t$ is general. 
Hence the vertical map on the right hand side 
is also an isomorphism by the diagram. 
Therefore, for general $t$, 
$h^0(i_*F\otimes_{{\cal O}_{\cal Z}}{\cal O}_{Z_t}) 
= h^0(i_*F_t)$.

 Let $\nu_t : {\tilde Z}_t \to Z_t$ be a 
resolution of singularities. 
By Theorem 4 
$i_*F_t \cong {\nu_t}_*\Omega^2_{{\tilde Z}_t}$. 
Since $h^0(i_*F_t) = h^0({\tilde Z}_t, 
\Omega^2_{{\tilde Z}_t}) = 
h^2({\tilde Z}_t, {\cal O}_{{\tilde Z}_t})  
= h^2(Z_t, {\cal O}_{Z_t})$, $h^0(i_*F_t)$ 
is a constant function of $t$. 
$h^0(i_*F\otimes_{{\cal O}_{\cal Z}}{\cal O}_{Z_t})$ 
is an upper semi-continous function of $t$. 
For any $t$, 
$h^0(i_*F\otimes_{{\cal O}_{\cal Z}}{\cal O}_{Z_t})
 \leq h^0(i_*F_t)$ and the 
equality holds for general $t$. 
Since $h^0(i_*F_t)$ is constant, this implies that 
$h^0(i_*F\otimes_{{\cal O}_{\cal Z}}{\cal O}_{Z_t})$ 
is constant and 
$h^0(i_*F\otimes_{{\cal O}_{\cal Z}}{\cal O}_{Z_t}) 
= h^0(i_*F_t)$ for all $t$. 
By a theorem of Grauert (cf. [Ha, Corollary 12.9]) 
$g_*(i_*F)$ is a locally free sheaf on $\Delta$ and 
the natural map 
$g_*(i_*F)\otimes_{{\cal O}_{\Delta}}k(0) 
\to H^0(i_*F\otimes_{{\cal O}_{\cal Z}}{\cal O}_{Z}) 
\cong H^0(i_*F_0)$ is an isomorphism. 
This implies that there is a lift $\tilde \omega 
\in \Gamma(g^{-1}(\Delta_{\epsilon}), 
\Omega^2_{{\cal U}/\Delta})$ of $\omega$. 
Let us consider $\wedge^{n/2}\omega$ 
as a section of the dualizing sheaf $\omega_Z$. 
Then $\wedge^{n/2}\tilde \omega$ can be regarded 
as a section of 
$\omega_{{\cal Z}_{\epsilon}/\Delta_{\epsilon}}$. 
Since $\wedge^{n/2}\omega$ generates the line bundle 
$\omega_Z$, $\wedge^{n/2}\tilde\omega$ also generates 
$\omega_{{\cal Z}_{\epsilon}/\Delta_{\epsilon}}$, 
if necessary, by taking $\epsilon$ smaller. 
Therefore $\omega_t := \tilde\omega\vert_{{\cal U}_t}$ 
is a non-degenerate 2-form for 
$t \in \Delta_{\epsilon}$. 
Since ${\cal Z}_t$ has only rational Gorenstein 
singularities, 
$(Z_t, \omega_t)$ is a symplectic variety by Theorem 6. 
\vspace{0.2cm}

{\bf Remark}. By virtue of Proposition 5, 
Theorem 7 seems true
 even when $g$ is a proper flat morphism and $Z$
 is a symplectic variety. 
The missing ingredients consist of two parts;  
(a) for a compact Kaehler space $Z$ 
with rational singularities, 
are the natural maps 
$H^i(Z, {\bold C}) \to H^i(Z, {\cal O}_Z)$ 
surjective for all $i$ ? 
(b) does Theorem 4 hold for an arbitrary 
rational Gorenstein singularity ? 
(in the proof of Theorem 4 we used a vanishing theorem 
of [St] and this vanishing theorem is only known 
for the case $X$ is embedded as an open subset 
in a complex projective variety.) 

If these two questions are affirmative, 
Theorem 7 holds in this full generality. \vspace{0.15cm}

The following will be used later. \vspace{0.15cm}

{\bf Theorem 7'}. {\em Let $(Z, \omega)$ be a 
symplectic variety with 
$\mathrm{codim}(\Sigma \subset Z) \geq 4$.  
Let $g : {\cal Z} \to \Delta^n$ be a proper flat 
morphism from 
${\cal Z}$ to a n-dimensional unit disc 
$\Delta^n$ with $g^{-1}(0) = Z$. 
Then $\omega$ extends sideways in the flat family 
so that it gives a symplectic 2-form 
$\omega_t$ on each fiber $Z_t$ for 
$t \in \Delta^n_{\epsilon}$ with a sufficiently small 
$\epsilon$.}  \vspace{0.15cm}

{\em Proof}. As in the proof of Theorem 7, we put 
${\cal U} := \{z \in {\cal Z}$; 
$g$ is smooth at $z$ $\}$. 
Denote by $i$ the natural inclusion of ${\cal U}$ to 
${\cal Z}$. 
Write $\pi$ for $g\circ i$. 
Put $F := \Omega^2_{{\cal U}/\Delta}$, $F_0 
:= \Omega^2_U$ and $U := {\cal U} \cap Z$. 
When $n > 1$, 
$i_*F$ is not necessarily flat over $\Delta^n$. 
Instead of using base change theorem 
we shall apply the comparison theorem between 
formal and analytic higher direct images by 
Banica [B-S, VI Proposition 4.2]. 

Let $t_1$, $t_2$, ..., $t_n$ be 
coordinates of $\Delta^n$. We assume 
that ${\cal Z}_p$ is Kaehler and 
$\mathrm{codim}({\mathrm{Sing}}({\cal Z}_p) 
\subset {\cal Z}_p) \geq 4$ for all 
point $p \in \Delta^n$. \vspace{0.12cm} 

{\bf (i)}:  For n-tuple of positive integers 
$(i_1, i_2, ..., i_n)$ 
, and for a point $p = (p_1, ..., p_n) \in \Delta^n$, 
we put 
$A_{(i_1, ..., i_n ; p)} = 
{\bold C}[t_1, ..., t_n]/
((t_1 - p_1)^{i_1}, (t_2 - p_2)^{i_2}, ..., 
(t_n - p_n)^{i_n})$ and 
$\Delta_{(i_1, ..., i_n ; p)} := 
\mathrm{Spec}A_{(i_1, ..., i_n ; p)}$. 
Define $F_{(i_1, ..., i_n ; p)} 
:= F \otimes_{{\cal O}_{\Delta^n}}{\cal O}_
{\Delta_{(i_1, ..., i_n ; p)}}$ and 
$U_{(i_1, ..., i_n ; p)} := 
{\cal U} \times_{\Delta^n}\Delta_{(i_1, ..., i_n ; p)}$. 
When $p = (0,0, ..., 0)$, we write $A_{(i_1, ..., i_n)}$ 
(resp. $\Delta_{(i_1, ..., i_n)}$, $U_{(i_1, ..., i_n)}$) 
for $A_{(i_1, ..., i_n ; p)}$ (resp. 
$\Delta_{(i_1, ..., i_n ; p)}$, $U_{(i_1, ..., i_n ; p)}$). 
\vspace{0.12cm}

{\bf (ii)}: For the later use we shall extend 
the notation above 
to the case where some indices $i_l$ are infinity. 
For simplicity, we assume that $i_1 = \infty$, ..., 
$i_{k-1} = \infty$ and $i_k$, ..., $i_n$ are positive 
integers. The notation in the general case would be 
clear from the explanation below.  
Write $\Delta^n = \Delta_{<k} \times \Delta_{\geq k}$, 
where $\Delta_{<k}$ is the $k-1$ dimensional polydisc 
with coordinates $t_1$, ..., $t_{k-1}$ and 
$\Delta_{\geq k}$ is the $n-k+1$ dimensional polydisc 
with coordinates $t_k$, ..., $t_n$.
For $n-k+1$-tuple of positive integers $(i_k, ..., i_n)$, 
and for a point $p = (p_{i_k}, ..., p_n) 
\in \Delta_{\geq k}$,  
we define 
$\Delta_{(\infty, ..., 
\infty, i_k, ..., i_n ;p)}$ as 
$\Delta_{<k} \times 
{\mathrm{Spec}}{\bold C}[t_k, ..., t_n]/
((t_k - p_k)^{i_k}, (t_{k+1} - p_{k+1})^{i_{k+1}}, ..., 
(t_n - p_n)^{i_n})$.  
Now  $F_{(\infty, ..., \infty, i_k, ..., i_n ; p)}$ and 
$U_{(\infty, ..., \infty, i_k, ..., i_n ; p)}$ are 
defined in a similar way as the case where 
$i_1, ..., i_n$ are all 
finite.  
We denote by ${\cal O}_{(i_1, ..., i_n ; p)}$ the 
structure sheaf of $\Delta_{(i_1, ..., i_n ; p)}$. 
\vspace{0.25cm}

{\bf Claim}. (1) {\em $R^j\pi_*F_{(i_1, ..., i_n ; p)}$ are 
coherent for $j = 0,1$ and for all 
$(i_1, ..., i_n)$ with $0 < i_k \leq \infty$  
$(1 \le k \leq n)$.} \vspace{0.15cm}

(2) {\em The natural maps 
$\pi_*F_{(i_1 +1, i_2, ..., i_n ; p)} \to 
\pi_*F_{(i_1, i_2, ..., i_n ; p)}$, 
$\pi_*F_{(i_1, i_2 + 1, ..., i_n ; p)} 
\to \pi_*F_{(i_1, i_2, ..., i_n ; p)}$, ..., 
$\pi_*F_{(i_1, i_2, ..., i_n + 1 ; p)} 
\to \pi_*F_{(i_1, i_2, ..., i_n ; p)}$ 
are surjective for all $(i_1, ..., i_n)$ 
with $0 < i_k < \infty$  $(1 \le k \leq n)$.} 
\vspace{0.15cm} 

(3) {\em $\pi_*F$ is locally free and 
$\pi_*F\otimes k(p) \cong H^0(U, F_p)$ for 
each $p \in \Delta^n$.} 
\vspace{0.15cm}

{\em Proof}. We shall first prove (1) and (2). 
We finally conclude (3) by combining the 
comparison theorem [B-S, VI, Proposition 4.2] 
with (1) and (2). 

({\bf 1}): Since $\mathrm{codim}(\Sigma \subset Z) 
\geq 3$ and $\mathrm{depth}(F_{(i_1, ..., i_n ; p)})_{,q} 
= \dim U_{(i_1, ..., i_n ; p)}$ for 
$q \in U_{(i_1, ..., i_n ; p)}$, 
$i_*F_{(i_1, ..., i_n ; p)}$ and 
$R^1i_*F_{(i_1, ..., i_n ; p)}$ are both coherent. 
By the exact sequence 

$$  0 \to R^1g_*(i_*F_{(i_1, ..., i_n ; p)}) 
\to R^1\pi_*F_{(i_1, ..., i_n ; p)} 
\to g_*(R^1i_*F_{(i_1, ..., i_n ; p)} 
\to R^2g_*(i_*F_{(i_1, ..., i_n ; p)}) $$ 
we know that $R^j\pi_*F_{(i_1, ..., i_n ; p)}$ 
are coherent for $j = 0,1$.  
\vspace{0.15cm}

({\bf 2}): We assume that $p = (0,0, ..., 0)$ 
because the proof are the same for all points $p$.  
Since $\mathrm{codim}(\Sigma \subset Z) 
\geq 4$, the spectral sequence 

$$  E_1^{p,q} := H^q(U, \Omega^p_U) 
=> H^{p+q}(U, {\bold C}) $$

degenerates at $E_1$ terms when $p + q = 2$ 
([Oh, Na 1, Lemma 2.5]). 
Let $U_m \to S_m$ be a flat deformation of $U$ over 
$S_m := \mathrm{Spec} A_m$ where 
$ A_m := {\bold C}[t]/(t^{m + 1})$ with 
$m \in {\bold Z}_{>0}$. Then, by [Na 1, Lemma 2.6], 
we see that the spectral sequence 

$$  E_1^{p,q} := H^q(U, \Omega^p_{U_m/S_m}) 
=> H^{p+q}(U, A_m) $$ 
degenerates at $E_1$ terms with $p + q = 2$. 
Here $A_m$ means the constant sheaf on $U$ with values in 
$A_m$. If we put $U_{m-1} := U_m \times_{S_m} S_{m-1}$, 
then the restriction map $H^q(U, \Omega^p_{U_m/S_m}) 
\to H^q(U, \Omega^p_{U_{m-1}/S_{m-1}})$ is surjective for 
$(p, q)$ with $p + q = 2$ ([Na 1, lemma 2.6]).     

To prove (2) we only have to check that 
$\pi_*F_{(i_1 +1, i_2, ..., i_n)} \to 
\pi_*F_{(i_1, i_2, ..., i_n)}$ is surjective by symmetry. 
One can split up the surjection 
$A_{(i_1 + 1, i_2, ..., i_n)} \to A_{(i_1, ..., i_n)}$ 
into a finite sequence of small extensions: 
$A_{(i_1 + 1, ..., i_n)} = A^{(N)} \to A^{(N-1)} 
\to ... A^{(1)} \to A_{(i_1, ..., i_n)}$ 
where $N = \dim_{\bold C} A_{(1, ..., i_n)}$. 
By definition, $K_{j} := \mathrm{ker}[A^{(j+1)} 
\to A^{(j)}]$ 
are one dimensional ${\bold C}$ vector spaces. 
For each $A^{(j+1)} \to A^{(j)}$ 
we can choose homomorphisms of local 
${\bold C}$ algebras  
${\bold C} [t]/(t^{m_j + 1}) \to A^{(j+1)}$ and 
${\bold C} [t]/(t^{m_j}) \to A^{(j)}$ in such a 
way that the diagram    

\begin{equation}
\begin{CD}
0 @>>> K_j @>>> A^{(j+1)} @>>> A^{(j)} @>>> 0 \\
@. @V{\phi_j}VV @VVV @VVV  \\
0 @>>> (t^{m_j}) @>>> {\bold C} [t]/(t^{m_j + 1}) 
@>>> {\bold C} [t]/(t^{m_j}) 
@>>> 0  
\end{CD}
\end{equation}
commutes and $\phi_j$ is an isomorphism. 
We put $Z_{(j)} := {\cal Z} \times_{\Delta^n} 
\mathrm{Spec} A^{(j)}$, $F_{(j)} := 
F \otimes_{{\cal O}_{\cal Z}}{\cal O}_{Z_{(j)}}$, 
$Z_{m_j} := {\cal Z} \times_{\Delta^n} \mathrm{Spec} 
{\bold C} [t]/(t^{m_j + 1})$ and 
$F_{m_j} := F 
\otimes_{{\cal O}_{\cal Z}}{\cal O}_{Z_{m_j}}$. 
By the previous observation, $\pi_*F_{m_j} 
\to \pi_*F_{m_j - 1}$ is surjective. 
We see that $\pi_*F_{(j+1)} \to \pi_*F_{(j)}$ 
is surjective by the commutative diagram. 
Hence we know that $\pi_*F_{(i_1 +1, i_2, ..., i_n)} 
\to \pi_*F_{(i_1, i_2, ..., i_n)}$ is surjective. 
\vspace{0.15cm}

({\bf 3}): We shall prove, by induction on $k$, that 
$\pi_*F_{(\infty, ..., \infty, i_k, ..., i_n ; p)}$ 
are free 
${\cal O}_{(\infty, ..., \infty, i_1, ..., i_n ; p)}$ 
modules and 
$\pi_*F_{(\infty, ..., \infty, i_k + 1, ..., i_n ; p)} 
\to \pi_*F_{(\infty, ..., \infty, i_k, ..., i_n ; p)}$ 
are surjective for 
all $n-k+1$ tuple 
$(i_k, ..., i_n)$ without infinity. For $k = 1$, 
they are nothing but (2) of Claim. 
Let 
$(\pi_*F_{(\infty, ..., \infty, i_k, ..., i_n ; p)})
\hat{ }_{p_{k-1}}$ 
be the completion of 
$\pi_*F_{(\infty, ..., \infty, i_k, ..., i_n ; p)}$ 
along the divisor $\{t_{k-1} = p_{k-1} \}$ of $\Delta^n$.

 It suffices to prove that 
$(\pi_*F_{(\infty, ..., 
\infty, i_k + 1, ..., i_n ; p)})
\hat{}_{p_{k-1}} 
\to (\pi_*F_{(\infty, ..., 
\infty, i_k, ..., i_n ; p)})
\hat{}_{p_{k-1}}$ 
are surjective for all $p_{k-1}$ 
in order to prove that 
$\pi_*F_{(\infty, ..., \infty, i_k + 1, ..., i_n ; p)} 
\to \pi_*F_{(\infty, ..., \infty, i_k, ..., i_n ; p)}$ 
are surjective. 
By the comparison theorem [B-S, VI, Proposition 4.2] 
and by (1), we have 

$$  (\pi_*F_{(\infty, ..., 
\infty, i_k, ..., i_n)})\hat{}_{p_{k-1}} 
\cong \underleftarrow{\lim} 
\pi_*F_{(\infty, ..., \infty, m, i_k, ..., i_n ; p_{k-1}, p)}. $$ 

Note that $p = (p_k, ..., p_n) \in  
\Delta_{\geq k}$, and $p_{k-1} \in \Delta^1(t_{k-1})$, 
where $\Delta^1(t_{k-1})$ is the 1-dimensional disc 
with a coordinate $t_{k-1}$. 
By the induction hypothesis, the maps 

$$\pi_*F_{(\infty, ..., \infty, m, i_k + 1, ..., i_n ; p_{k-1}, p)} 
\to \pi_*F_{(\infty, ..., \infty, m, i_k, ..., i_n ; p_{k-1}, p)}$$ 
and 
$$\pi_*F_{(\infty, ..., \infty, m +1, 1, ..., i_n ; p_{k-1}, p)} 
\to \pi_*F_{(\infty, ..., \infty, m, 1, ..., i_n ; p_{k-1}, p)} $$ 
are both surjective for all $p_{k-1}$. 
Therefore we conclude that 

$$(\pi_*F_{(\infty, ..., \infty, i_k + 1, ..., i_n)})\hat{}
_{p_{k-1}} 
 \to 
(\pi_*F_{(\infty, ..., \infty, i_k, ..., i_n)})\hat{}
_{p_{k-1}}$$ 
are surjective for all $p_{k-1}$. \vspace{0.12cm}
 
We shall next prove that 
$\pi_*F_{(i_1, ..., i_n ; p)}$ are free 
${\cal O}_{(i_1, ..., i_n ; p)}$ 
module. 

We shall prove that 
$\pi_*F_{(\infty, ..., \infty, i_k, ..., i_n ; p)}$ 
is a free ${\cal O}_{(\infty, ..., \infty, i_k, ..., i_n ; p)}$ 
module by assuming that 
$\pi_*F_{(\infty, ..., \infty, *_{k-1}, ..., *_{n}; p')}$ 
are free 
${\cal O}_{(\infty, ..., \infty, *_{k-1}, ..., *_n ; p')}$ 
modules for all $*_j \in {\bold Z}_{>0}$ and all $p' \in 
\Delta_{\geq k-1}$.

We use the induction on the lexicographic order of 
$(i_k, ..., i_n)$. 
 First $\pi_*F_{(\infty, ..., \infty, 1, ..., 1 ; p)}$ 
is a free ${\cal O}_{(\infty, ..., \infty, 1, ..., 1 ; p)}$ module; 
in fact, let 
$(\pi_*F_{(\infty, ..., \infty, 1, ..., 1; p)})
\hat{ }_{p_{k-1}}$ 
be the completion of 
$\pi_*F_{(\infty, ..., \infty, 1, ..., 1;p)}$ 
along the divisor $\{t_{k-1} = p_{k-1} \}$ of $\Delta^n$. Then 

$$  (\pi_*F_{(\infty, ..., \infty, 1, ..., 1; p)})\hat{ }
_{p_{k-1}} 
\cong \underleftarrow{\lim} 
\pi_*F_{(\infty, ..., \infty, m, 1, ..., 1; p_{k-1}, p)}. $$

Since $\pi_*F_{(\infty, ..., \infty, m, 1, ..., 1; p_{k-1}, p)}$ 
are free 
${\cal O}_{(\infty, ..., \infty, m, 1, ..., 1; p_{k-1}, p)}$ modules
 by assumption and 
$\pi_*F_{(\infty, ..., \infty, m, 1, ..., 1; p_{k-1}, p)} 
\to \pi_*F_{(\infty, ..., \infty, m - 1, 1, ..., 1; p_{k-1}, p)}$ 
are surjective, the right hand side 
is a free 
$\hat{\cal O}_{(\infty, ..., \infty, 1, ..., 1; p)}$ module. 
Therefore, $\pi_*F_{(\infty, ..., \infty, 1, ..., 1; p)}$ 
is a free 
${\cal O}_{(\infty, ..., \infty, 1, ..., 1; p)}$ module. 

Next consider 
$\pi_*F_{(\infty, ..., \infty, i_k, ..., i_n; p)}$. 
By induction 
$\pi_*F_{(\infty, ..., \infty, i_k - 1, ..., i_n; p)} = 
\pi_*F_{(\infty, ..., \infty, i_k, ..., i_n; p)}
\otimes {\cal O}_{\infty, ..., \infty, i_k - 1, ..., i_n; p)}$ 
is isomorphic to 
$({\cal O}_{\infty, ..., \infty, i_k - 1, ..., i_n; p)})^r$. 
By Nakayama's lemma, we have a surjection from 
$({\cal O}_{\infty, ..., \infty, i_k, ..., i_n; p)})^r$ to 
$\pi_*F_{(\infty, ..., \infty, i_k, ..., i_n; p)}$. 
We then have a commutative diagram with exact rows:

\begin{equation}
\begin{CD}
 ({\cal O}_{(..., \infty, 1, i_{k+1}, ..., i_n; p)})^r 
@>>> ({\cal O}_{(..., \infty, i_k, i_{k+1}, ..., i_n; p)})^r 
@>>> 
({\cal O}_{(..., \infty, i_k - 1, i_{k+1}, ..., i_n; p)})^r \\
@VVV @VVV @VVV   \\
 \pi_*F_{(..., \infty, 1, ..., i_n; p)} @>>> 
\pi_*F_{(..., \infty, i_k, ..., i_n; p)} @>>> 
\pi_*F_{(..., \infty, i_k - 1, i_{k+1}, ..., i_n; p)} 
\end{CD}
\end{equation}

On each row, the first map is injective and 
the second one is surjective.
The right vertical map is an isomorphism. 
The middle vertical map is surjective; 
hence, by the Snake Lemma, 
the left vertical map is surjective. 
On the other hand, by the induction hypothesis, 
$\pi_*F_{(..., \infty, 1, ..., i_n; p)}$ is a free 
${\cal O}_{(..., \infty, 1, ..., i_n; p)}$ module of rank $r$. 
This implies that the 
left vertical map is an isomorphism. 
Again, by the Snake lemma, the middle vertical map 
is an injection, hence an isomorphism.   Q.E.D.  
\vspace{0.2cm}

{\em Proof of Theorem 7' continued}. By Claim (3) the 
symplectic 2-form $\omega$  on $U$ extends sideways. 
Since Theorem 4 (hence Theorem 6) holds for a 
(non-algebraic) singularity with 
$\mathrm{codim}(\Sigma \subset X) \geq 4$ by [Fl], 
the rest of the argument is the same as Theorem 7. 
Q.E.D. \vspace{0.2cm}

We now consider the following situation:  
Let $Z$ be a symplectic variety.  Put 
$\Sigma := \mathrm{Sing}(Z)$ and 
$U := Z \setminus \Sigma$.  
Let $\overline \pi : {\cal Z} \to S$ be 
the Kuranishi family of $Z$, which is, by definition, 
a semi-universal flat deformation of $Z$ with 
$\overline \pi^{-1}(0) = Z$ for the reference point 
$0 \in S$. When 
$\mathrm{codim}(\Sigma \subset Z) \geq 4$, 
$S$ is smooth by [Na 1, Theorem 2.4]. 
Define ${\cal U}$ to be the locus in ${\cal Z}$ where 
$\overline \pi$ is a smooth map and let 
$\pi : {\cal U} \to S$ be the restriction of 
$\overline \pi$ to ${\cal U}$. 
The following is a generalization of 
the Local Torelli Theorem [Be 2, Theoreme 5] 
to singular symplectic varieties.
\vspace{0.2cm}

{\bf Theorem 8} {\em Assume that $Z$ is a 
${\bold Q}$-factorial projective symplectic variety. 
Assume $h^1(Z, {\cal O}_Z) = 0$, 
$h^0(U, \Omega^2_U) = 1$, $\dim Z = 2l \geq 4$ and 
$\mathrm{Codim}(\Sigma \subset Z) \geq 4$. 
Then the following hold.} \vspace{0.15cm}

(1) {\em $R^2\pi_*(\pi^{-1}{\cal O}_S)$ is a free 
${\cal O}_S$ module of finite rank. 
Let ${\cal H}$ be the image of the composite 
$R^2{\overline \pi}_*{\bold C} \to R^2\pi_*{\bold C} 
\to R^2\pi_*(\pi^{-1}{\cal O}_S)$. 
Then ${\cal H}$ is a local system on $S$ with 
${\cal H}_s = 
H^2({\cal U}_s, {\bold C})$ for $s \in S$.} 
\vspace{0.15cm}

(2) {\em The restriction map $H^2(Z, {\bold C}) 
\to H^2(U, {\bold C})$ is an isomorphism. 
Take a resolution $\nu: \tilde Z \to Z$ 
in such a way that $\nu^{-1}(U) \cong U$. 
For $\alpha \in H^2(U, {\bold C})$ we take a lift 
$\tilde\alpha \in H^2(\tilde Z, {\bold C})$ 
by the composite $H^2(U, {\bold C}) 
\cong H^2(Z, {\bold C}) \to H^2(\tilde Z, {\bold C})$. 
Choose $\omega \in H^0(U, \Omega^2_U) = {\bold C}$. 
This $\omega$ extends to a holomorphic 2-form on 
$\tilde Z$. Normalize $\omega$ in such a way that 
${\int_{\tilde Z} (\omega \overline \omega)^l} = 1$.  
Then one can define a quadratic form 
$q : H^2(U, {\bold C}) \to {\bold C}$ as} 

$$  q(\alpha) :=  l/2 {\int_{\tilde Z}
(\omega \overline \omega)^{l-1}\tilde\alpha^2} + 
(1-l)(\int_{\tilde Z}\omega^l
\overline\omega^{l-1}\tilde\alpha)
(\int_{\tilde Z}\omega^{l-1}
\overline\omega^l\tilde\alpha ). $$

{\em This form is independent of the choice of 
$\nu : \tilde Z \to Z$.} \vspace{0.15cm}

(3) {\em Put $H := H^2(U, {\bold C})$. 
Then there exists a trivialization of the local system 
${\cal H}$:  ${\cal H} \cong H \times S$. 
Let 
${\cal D} := \{ x \in {\bold P} (H) ; 
q(x) =0, q(x + \overline x) > 0 \}$. 
Then one has a period map 
$p: S \to {\cal D}$ and $p$ is a local isomorphism.} 
\vspace{0.2cm}

{\bf Remark}. (1) The assumption that 
$\mathrm{Codim}(\Sigma \subset Z) \geq 4$ 
is always satisfied when $Z$ has only 
terminal singularities [Na 2, Theorem]. 

(2) One can apply Theorem 8 for $Z$ for 
irreducible symplectic V-manifolds (cf. [Fu]); but 
the results seems rather trivial by Schlessinger's 
rigidity theorem for quotient singularities. 

Most interesting objects are 
$\overline{M_{0,c}}$ in Example (iii) 
($c \geq 6$) of Introduction. $\overline{M_{0,c}}$
conjecturelly satisfy the assumption of Theorem 8. 
Since they have more complicated singularities than 
quotient singularities, one expects that they give 
non-trivial examples for Theorem 8. Actually Theorem 8 
is motivated by them. \vspace{0.15cm}
      
{\em Proof of (1)}: The first statement follows 
from the following result. \vspace{0.15cm}

{\bf Proposition 9}. {\em The spectral sequence}

$$ E^{p,q}_1 = R^q\pi_*\Omega^p_{{\cal U}/S} 
=> R^{p+q}\pi_*(\pi^{-1}{\cal O}_S) $$

{\em degenerates at $E_1$ terms for $p + q = 2$. 
Moreover, $E^{p,q}_1$ are locally free sheaves 
for $p + q = 2$.} \vspace{0.15cm}

{\em Proof of Proposition 9}. (a): $E^{2,0}_1$ 
is locally free and compatible with base change 
by Theorem 7'. Since $Z$ is Gorenstein and 
$\mathrm{Codim}(\Sigma \subset Z) \geq 4$, 
$H^2({\cal Z}_s, {\cal O}_{{\cal Z}_s}) 
\cong H^2({\cal U}_s, {\cal O}_{{\cal U}_s})$. 
Therefore $E^{0,2}_1$ is also locally free and 
compatible with base change by the proof of Theorem 7.  
By [Na 1, Theorem 2.4] $S$ is amooth. 
Since $h^1(Z, {\cal O}_Z) = 0$, 
we have $h^0(Z, \Theta_Z) = 0$, hence 
$\pi: {\cal Z} \to S$ is the universal family. 
This implies that $T_{S,s} \cong 
H^1({\cal U}_s, \Theta_{{\cal U}_s})$ for $s \in S$, 
where $T_{S,s}$ is the tangent space of $S$ at $s$. 
On the other hand, there are natural identifications 
$H^1({\cal U}_s, \Theta_{{\cal U}_s}) 
\cong H^1({\cal U}_s, \Omega^1_{{\cal U}_s})$ 
by a relative symplectic 2-form of $\pi$ 
(such a 2-form exists by Theorem 7'). Therefore, 
$E^{1,1}_1$ is locally free 
and compatible with base change. 

(b): We shall prove that the composed map 
$H^2({\cal Z}, {\bold C}) \to H^2({\cal U}, {\bold C}) 
\to H^2(U, {\bold C})$ is surjective when 
we choose $S$ small enough. 
Since $H^2({\cal Z}, {\bold C}) 
\cong H^2(Z, {\bold C})$, 
it suffices to show that the restriction map 
$\alpha: H^2(Z, {\bold C}) \to H^2(U, {\bold C})$ 
is surjective. 
Let $\nu: \tilde Z \to Z$ be a resolution of 
singularities whose exceptional locus $E$ 
is a simple normal crossing divisor. 
Then $\alpha$ can be factorized as 
$H^2(Z, {\bold C}) \to H^2(\tilde Z, {\bold C}) 
\to H^2(U, {\bold C})$. 
There is an exact sequence of mixed Hodge structures 

$$  H^2_E(\tilde Z, {\bold C}) \to 
H^2(\tilde Z, {\bold C}) \to 
H^2(U, {\bold C}) \to H^3_E(\tilde Z, {\bold C}). $$ 

Note that $H^3_E(\tilde Z)$ has the 
mixed Hodge structure 
with weights $\geq 3$. 
On the other hand, $H^2(U)$ has the pure 
Hodge structure \footnote
{Since $\mathrm{Codim}(\Sigma \subset Z) \geq 4$, 
by [Oh, Na 1, Lemma 2.5] the Hodge spectral sequence 
$H^q(U, \Omega^p_U) => H^{p+q}(U, {\bold C})$ 
degenerates at $E_1$ terms when $p + q = 2$. 
Moreover, $h^2(U, {\cal O}_U) = h^0(U, \Omega^2_U)$. 
This filtration coincides with the Hodge filtration 
of the mixed Hodge structure on $H^2(U, {\bold C})$. 
In fact, there is a natural map $\phi_{p,q} : 
H^q({\tilde Z}, \Omega^p_{\tilde Z}(\log E)) 
\to H^q(U, \Omega^p_U)$ 
for each $p$ and $q$. 
By [Fl] or Proposition 1, we have $\nu_*\Omega^2_
{\tilde Z}(\log E) \cong i_*\Omega^2_U$, hence 
$\phi_{2,0}$ is an isomorphism
. On the other hand, since $Z$ has only 
rational singularities and 
$\mathrm{Codim}(\Sigma \subset Z) \geq 4$, 
$\phi_{0,2}$ is an isomorphism. 
Recall that the spectral sequence 
$H^q({\tilde Z}, \Omega^p_{\tilde Z}(\log E)) 
=> H^{p+q}({\tilde Z}, \Omega^{\cdot}_{\tilde Z}
(\log E)) = H^{p+q}(U, {\bold C})$ degenerates 
at $E_1$ terms. Then we 
know that $\phi_{1,1}$ is also an isomorphism. 
Thus our filtration coincides with the Hodge 
filtration of the mixed Hodge structure of 
$H^2(U, {\bold C})$. Since 
$h^2(U, {\cal O}_U) = h^0(U, \Omega^2_U)$, 
the mixed Hodge structure is pure.}    
of weight 2 because 
$\mathrm{Codim}(\Sigma \subset Z) \geq 4$. 
Therefore 
$H^2(\tilde Z, {\bold C}) \to H^2(U, {\bold C})$ 
is a surjection. 

  We shall prove that, as a ${\bold C}$ vector space, 
$H^2(\tilde Z, {\bold C})$ is generated by the image of 
$H^2(Z, {\bold C})$ and the image of 
$H^2_E(\tilde Z, {\bold C})$; if this is proved, 
then $\alpha$ is surjective by the exact sequence.  
Let $E = \Sigma E_i$ be the irreducible decomposition 
and denote by $[E_i] \in H^2(\tilde Z, {\bold C})$ 
the cohomology class corresponding to 
the divisor $E_i$. 
Then we have $H^2_E(\tilde Z, {\bold C}) 
\cong \oplus {\bold C}[E_i]$. 
Since $Z$ is ${\bold Q}$-factorial and $Z$ has 
only rational singularities, 
$\mathrm{im}[H^2(\tilde Z, {\bold C}) 
\to H^0(Z, R^2\nu_*{\bold C})] = 
\mathrm{im}[\oplus {\bold C}[E_i] 
\to H^0(Z, R^2\nu_*{\bold C})]$ by [Ko-Mo, (12.1.6)]. 
Note that $H^2(Z, {\bold C}) = 
\mathrm{Ker}[H^2(\tilde Z, {\bold C}) 
\to H^0(Z, R^2\nu_*{\bold C})]$ because 
$Z$ has only rational singularities. 
The argument here shows that the restriction map 
$H^2(Z, {\bold C}) \to H^2(U, {\bold C})$ is, in fact, 
an isomorphism because $H^2(Z, {\bold C}) 
\to H^2(\tilde Z, {\bold C})$ is an injection. 
(Note that $R^1\nu_*{\bold C} = 0$ because 
$Z$ has rational singularities.)
         
 (c): Let $k(0)$ be the skyscraper sheaf 
supported at $0 \in S$ which is defined as the quotient 
${\cal O}_S/m_{0}$, where $m_0$ is the ideal sheaf of 
$0 \in S$. 
Then the natural map 
$(R^2\pi_*(\pi^{-1}{\cal O}_S))_0 
\to R^2\pi_*(\pi^{-1}k(0))$ is surjective. 
In fact, $(R^2\pi_*(\pi^{-1}{\cal O}_S))_0$ 
factors the map $(R^2\overline\pi_*{\bold C})_0 
\to R^2\pi_*(\pi^{-1}k(0))$. 
This map is nothing but the map 
$H^2({\cal Z}, {\bold C}) \to H^2(U, {\bold C})$  
which is a surjection by (b). \vspace{0.2cm}

We are now in a position to prove the 
$E_1$ degeneracy of the spectral sequence. 
Let $0 \subset {\cal F}^2 \subset {\cal F}^1 
\subset {\cal F}^0 = R^2\pi_*(\pi^{-1}{\cal O}_S)$ 
be the decreasing filtration defined by 
the spectral sequence. By checking the coherence 
of each $E^{p,q}_k$ term 
(under the assumption 
$\mathrm{codim}(\Sigma \subset Z) \geq 4$), 
we can see that ${\cal F}^i$ are coherent. 
We shall prove that 
$(Gr^i_{\cal F})_0 = 
(R^{2-i}\pi_*\Omega^i_{{\cal U}/S})_0$ for 
$i = 0, 1, 2$. 

($ i = 0$): Since $E^{0,2}_{\infty} \subset E^{0,2}_1$, 
it is enough to prove that the map 
$(R^2\pi_*(\pi^{-1}{\cal O}_S))_0 
\to (R^2\pi_*{\cal O}_{\cal U})_0$ is surjective. 
We have a commutative diagram 

\begin{equation}
\begin{CD}
(R^2\pi_*(\pi^{-1}{\cal O}_S))_0 @>>> 
(R^2\pi_*{\cal O}_{\cal U})_0 \\
@VVV @VVV   \\
H^2(U, {\bold C}) @>>> H^2(U, {\cal O}_U) 
\end{CD}
\end{equation}
 
By (c) the vertical map on the left hand side 
is surjective. By (a), 
$(R^2\pi_*{\cal O}_{\cal U})_0\otimes k(0) 
\cong H^2(U, {\cal O}_U)$. 
Since $\mathrm{codim}(\Sigma \subset Z) \geq 4$, 
the spectral sequence 

$$ H^q(U, \Omega^p_U) => H^{p+q}(U, {\bold C}) $$ 
degenerates at $E_1$ terms when $p + q = 2$ 
([Oh, Na 1, Lemma 2.5]). 
Hence the horizontal map at the bottom is surjective. 
By Nakayama's Lemma we see that the horizontal map 
on the top is also surjective. \vspace{0.15cm}

($i = 1$): By the assumption, $E^{0,1}_1 = 0$, 
hence $E^{1,1}_{\infty} \subset E^{1,1}_1$. 
It is enough to prove that 
$({\cal F}^1)_0 \to (R^1\pi_*\Omega^1_{{\cal U}/S})_0$ 
is surjective. We have two commutative diagrams:

\begin{equation}
\begin{CD}
({\cal F}^1)\otimes k(0) @>>> 
(R^2\pi_*(\pi^{-1}{\cal O}_S)\otimes k(0) 
@>>> (R^2\pi_:{\cal O}_{\cal U})\otimes k(0) @>>> 0 \\
 @VVV @VVV @VVV  \\
F^1 @>>> H^2(U, {\bold C}) @>>> H^2(U, {\cal O}_U) @>>> 0 
\end{CD}
\end{equation}

\begin{equation}
\begin{CD}
({\cal F}^1)_0 @>>> (R^1\pi_*\Omega^1_{{\cal U}/S})_0 \\
@VVV @VVV   \\
F^1 @>>> H^1(U, \Omega^1_U). 
\end{CD}
\end{equation}
   
The $F^1$ in the first diagram is the filtration of 
$H^2(U, {\bold C})$ induced by the spectral sequence 

$$ H^q(U, \Omega^p_U) => H^{p+q}(U, {\bold C}). $$

The rows in the first diagram are exact. Moreover, 
$F_1 \to H^1(U, \Omega^1_U)$ is injective 
by definition of $F_1$. 
Let us look at the first diagram. 
By (a) the vertical map on 
the right hand side is an isomorphism. 
The middle vertical map is surjective by (c). 
Hence the vertical map on 
the left hand side is surjective.  

We next observe the second diagram. 
The map $({\cal F}^1)_0 \to F^1$ is 
surjective because it is factorized as 
$({\cal F}^1)_0 \to ({\cal F}^1)\otimes k(0) \to F^1$. 
Since $Gr_F^1 = H^1(U, \Omega^1_U)$, the 
horizontal map at the bottom is surjective. Since 
$(R^1\pi_:\Omega^1_{{\cal U}/S})_0 \otimes k(0) 
\cong H^1(U, \Omega^1_U)$, 
the map $({\cal F}^1)_0 
\to (R^1\pi_*\Omega^1_{{\cal U}/S})_0$ is 
surjective by Nakayama's lemma. \vspace{0.15cm}

($i = 2$): By the assumption $E^{1,0}_1 = 0$ and 
$E^{0,1}_1 = 0$; hence 
$E^{2,0}_{\infty} \subset E^{2,0}_1$. 
We shall prove that $({\cal F}^2)_0 \to 
(\pi_*\Omega^2_{{\cal U}/S})_0$ is surjective. 
We have two commutative diagrams:

\begin{equation}
\begin{CD}
({\cal F}^2)\otimes k(0) @>>> ({\cal F}^1)\otimes k(0) 
@>>> (R^1\pi_*\Omega^1_{{\cal U}/S})\otimes k(0) 
@>>> 0 \\
 @VVV @VVV @VVV  \\
F^2 @>>> F^1 @>>> H^1(U, \Omega^1_U) @>>> 0 
\end{CD}
\end{equation}

\begin{equation}
\begin{CD}
({\cal F}^2)_0 @>>> (\pi_*\Omega^2_{{\cal U}/S})_0 \\
@VVV @VVV   \\
F^2 @>>> H^0(U, \Omega^2_U). 
\end{CD}
\end{equation}
   
In the first diagram the vertical map on 
the right hand side is an isomorphism, 
and the middle vertical map is surjective 
by the argument of ($i = 1$). 
Therefore ${\cal F}^2 \otimes k(0) \to F^2$ 
is surjective; this implies that, 
in the second diagram, 
the vertical map on the left hand side is surjective. 
Look at the second diagram. 
Since the horizontal map at the bottom is an 
isomorphism and the vertical map on the right hand side 
is surjective by (a), we conclude that 
$({\cal F}^2)_0 \to (\pi_*\Omega^2_{{\cal U}/S})_0$ 
is surjective by Nakayama's lemma. Q.E.D. 
\vspace{0.15cm}

Let us prove that ${\cal H}$ is a local system on $S$ 
with ${\cal H}_s = H^2({\cal H}_s, {\bold C})$ 
for $s \in S$. 
First note that $R^2\pi_*(\pi^{-1}{\cal O}_S)\otimes k(s) 
\cong H^2({\cal U}_s, {\bold C})$ by Lemma 9 because 
$R^q\pi_*\Omega^p_{{\cal U}/S}$, $p + q = 2$ are 
compatible with base change, and the spectral sequence 
$H^q({\cal U}_s, \Omega^p_{{\cal U}_s}) => 
H^2({\cal U}_s, {\bold C})$ 
degenerates at $E_1$ terms with 
$p + q = 2$. We take $S$ small enough so that 
${\cal Z}$ has strong deformation retract to $Z$. 
Choose $s \in S$. Then we have a diagram:

$$  H^2(U, {\bold C}) \leftarrow H^2(Z, {\bold C}) 
\rightarrow H^2({\cal Z}_s, {\bold C}) 
\rightarrow H^2({\cal U}_s, {\bold C}). $$    

The map $H^2(Z, {\bold C}) \to H^2(U, {\bold C})$ 
is an isomorphism by (b); hence we have a map 
$\phi_s : H^2(U, {\bold C}) 
\to H^2({\cal U}_s, {\bold C})$. 
Consider the map 
$\iota: \Gamma (S, R^2{\overline \pi}_*{\bold C}) 
\to \Gamma (S, R^2\pi_*(\pi^{-1}{\cal O}_S))$ 
induced by the sheaf homomorphism 
$R^2{\overline \pi}_*{\bold C} \to R^2\pi_*{\bold C} 
\to R^2\pi_*(\pi^{-1}{\cal O}_S)$. 
Denote by $\iota (s)$ the composite of 
$\iota$ with the evaluation map at $s$: 
$\Gamma (S, R^2\pi_*(\pi^{-1}{\cal O}_S)) 
\to H^2({\cal U}_s, {\bold C})$. 
For $\sigma \in 
\Gamma (S, R^2\pi_*(\pi^{-1}{\cal O}_S))$ 
we see that $(\iota (s))(\sigma) = 
\phi_s((\iota (0))(\sigma))$. 
If we choose $\sigma_1, ..., \sigma_r \in 
\Gamma(S, R^2{\overline \pi}_*{\bold C})$ in such a way 
that $(\iota (0))(\sigma_1), ..., (\iota (0))(\sigma_r)$ 
span $H^2(U, {\bold C})$ as a ${\bold C}$ vector space, 
then $(\iota (s))(\sigma_1), ..., (\iota (s))(\sigma_r)$ 
also span 
$H^2({\cal U}_s, {\bold C})$ for each $s \in S$, 
if necessary, by replacing $S$ by a smaller one. 
This implies that $H^2({\cal Z}_s, {\bold C}) 
\to H^2({\cal U}_s, {\bold C})$ 
is surjective in the diagram. 
Since $H^2({\cal Z}_s, {\bold C}) \to 
H^2({\cal U}_s, {\bold C})$ is injective 
(cf. the final statement of (b)), it is an isomorphism. 
Moreover, since $\dim H^2({\cal U}_s, {\bold C})$ are 
constant, we know by the diagram that 
$H^2(Z, {\bold C}) \cong H^2({\cal Z}_s, {\bold C})$. 
This completes the proof of (1) of Theorem 8. 
\vspace{0.15cm}

{\em Proof of (2) and (3)}: The statement of (2) is 
now clear from the arguments in the proof of (1). 
In the proof of (1), we have constructed 
an isomorphism $\phi_s : H^2(U, {\bold C}) 
\to H^2({\cal U}_s, {\bold C})$ for each 
$s \in S$. The trivialization $H \times S \to \cal{H}$ 
is given by $(x, s) \to \phi_s(x) \in 
{\cal H}_s = H^2({\cal U}_s, {\bold C})$. 
By Theorem 7' and the assumption, 
$\pi_*\Omega^2_{{\cal U}/S}$ is a line bundle on $S$ and 
is compatible with base change. 
Then we can choose $\tilde\omega \in 
\Gamma (S, \pi_*\Omega^2_{{\cal U}/S})$ 
in such a way that 
$\tilde\omega_0 = \omega$ and 
${\int_{{\cal U}_s} (\tilde{\omega}_s 
\overline{\tilde\omega}_s)^l} = 1$ for $s \in S$. 
By Hodge decomposition 
$H^2({\cal U}_s, {\bold C}) = 
H^0({\cal U}_s, \Omega^2_{{\cal U}_s}) 
\oplus H^1({\cal U}_s, \Omega^1_{{\cal U}_s}) 
\oplus H^2({\cal U}_s, {\cal O}_{{\cal U}_s})$, 
$\tilde\omega_s$ becomes an element of 
$H^2({\cal U}_s, {\bold C})$. 
Now the period map $p: S \to {\bold P}(H)$ is given by 
$s \to [\phi_s^{-1}(\tilde\omega_s)]$. 
The image of $p$ is contained in 
$\cal D$, and $p$ is a local isomorphism between 
$S$ and $\cal D$. The proofs are similar to 
[Be 2, Theoreme 5]. \vspace{0.15cm}

{\bf Remark (1)}:  In (2) the quadratic form 
$q(\alpha)$ can be defined by an arbitrary lift 
${\tilde \alpha} \in H^2(\tilde Z, {\bold C})$ of 
$\alpha \in H^2(U, {\bold C})$. 
The proof is as follows. 
Note that (cf. (b) of the proof of Proposition 9) 

$$ \mathrm{Ker}[H^2(\tilde Z, {\bold C}) 
\to H^2(U, {\bold C})] = 
\mathrm{im}[\oplus_i {\bold C} [E_i] \to  
H^2(\tilde Z, {\bold C})]. $$

Since 
$\int_{\tilde Z}\omega^l\overline\omega^{l-1}[E_i] 
= \int_{\tilde Z}\omega^{l-1}\overline\omega^l[E_i] = 0$, 
we only have to prove that $\int_{\tilde Z} 
(\omega \overline \omega)^{l-1}[E_i]\beta = 0$, 
$\beta \in H^2(\tilde Z, {\bold C})$. 
It is enough to prove that $\omega^{l-1}\vert_{E_i} = 0$. 
Set $S_i = \nu (E_i)$. We blow up $\tilde Z$ further 
and replace $\nu$ by a new resolution for which 
the inverse image of $S_i$ is a simple normal 
crossing divisor. 
Hereafter we call this new resolution 
$\nu: \tilde Z \to Z$ and put $F := \nu^{-1}(S_i)$. 
By definition, $F$ contains an irreducible component 
which is birational to the original $E_i$. 
By abuse of notation we call this component $E_i$. 
We only have to check that 
$\omega^{l-1}\vert_{E_i} = 0$ for the new $E_i$. 
We shall derive a contradiction by assuming that 
$\omega^{l-1}\vert_{E_i} \neq 0$.
Consider the map $F \to S_i$ induced by 
$\nu$. For $p \in S_i$, denote by 
$F_p$ the fiber over $p$. 
For a general point $p \in S_i$, 
$F_p$ is a normal crossing variety. 
Since $\mathrm{Codim}(S_i \subset Z) \geq 3$ 
by assumption, we have no non-zero holomorphic 
$2l - 2$ forms on $S_i$. 
Therefore, if $\omega^{l-1}\vert_{E_i} \ne 0$, 
then, for a general point $p \in S_i$, 
$H^0(F_p, \Omega^j_{F_p}) \neq 0$ for some $i > 0$.  
Put $k = \dim S_i$ and take a general complete 
intersection of $Z$ by $k$ hyperplanes: 
$H = H_1 \cap ... \cap H_k$. 
Put $\tilde H := \nu^{-1}(H)$. 
$H$ has canonical singularities, 
hence has rational singularities. 
Moreover, $f: \tilde H \to H$ is a 
resolution of singularities. 
Choose a point $p_i$ from $H \cap S_i$. 
We may assume that this $p_i \in S_i$ is 
general in the above sense. 
Note that $f^{-1}(p_i) = F_{p_i}$.
Since $R^jf_*{\cal O}_{\tilde H}= 0$ for $j > 0$, 
we see that $H^j(F_{p_i}, {\cal O}_{F_{p_i}}) = 0$ 
for $j > 0$. By the mixed Hodge structure on 
$H^j(F_{p_i})$ we conclude that 
$H^0(F_{p_i}, \hat\Omega^j_{F_{p_i}}) = 0$ for 
all $j > 0$. This is a contradiction. \vspace{0.15cm}

{\bf Remark (2)}: In Theorem 8, if we replace 
the ${\bold Q}$-factoriality condition by the 
next condition (*), then it is also valid for 
a non-projective symplectic variety: 

(*) $\mathrm{im}[H^2(\tilde Z, {\bold Q}) 
\to H^0(Z, R^2\nu_*{\bold Q})] = 
\mathrm{im}[\oplus {\bold Q}[E_i] 
\to H^0(Z, R^2\nu_*{\bold Q}].$

This condition is equivalent to the 
${\bold Q}$-factoriality when $Z$ is 
projective [Ko-Mo, (12.1.6)]. 
But when $Z$ is non-projective, 
they do not seem equivalent; for example, 
when $Z$ has no Weil divisors, 
${\bold Q}$-factoriality is 
meaningless. The condition (*) is 
an open condition for a family 
of symplectic varieties with 
terminal singularities. 
\vspace{0.2cm}

\begin{center}
Department of Mathematics, Graduate school of science, 
Osaka University, Toyonaka, Osaka 560, Japan 
\end{center}

\end{document}